\documentclass[aos]{imsart}

\RequirePackage{amsthm,amsmath,amssymb}
\RequirePackage[numbers]{natbib}
\usepackage{amsmath}
\usepackage{color, graphicx}
\usepackage[colorinlistoftodos]{todonotes}
\usepackage{xcolor}
\usepackage{verbatim}
\usepackage{natbib}
\usepackage{color,soul}
\usepackage{amsthm}
\usepackage{amsmath}
\usepackage{url} 
\usepackage{tikz}
\usepackage{etoolbox}
\usepackage{float}
\usepackage{amsmath}
\usepackage{float}
\usepackage{siunitx}
\usepackage{amssymb}
\usepackage{scalerel}
\usepackage{comment}
\usepackage{graphicx}
\usepackage[parfill]{parskip}
\usepackage{marvosym}
\usepackage{wrapfig}
\usepackage{resizegather}
\usepackage[toc,page]{appendix}
\usepackage{bbm}
\usepackage{graphicx}
\usepackage{amsbsy}
\usepackage{lipsum}
\usepackage{multicol}
\usepackage{amsthm}
\usepackage{pifont}
\usepackage{bbm}
\usepackage{bm}
\usepackage{dsfont}
\usepackage{bigints}
\usepackage{mathtools}
\usepackage{graphicx}
\usepackage{svg}
\usepackage{calc}
\usepackage{accents}
\usepackage{mwe}
\usepackage{graphicx,wrapfig,lipsum}
\usepackage{mathtools}
\usepackage{mathrsfs}
\usepackage{graphicx}
\usepackage[numbers]{natbib}
\usepackage[utf8]{inputenc}
\usepackage{enumitem}
\usepackage{letltxmacro}
\usepackage{nameref,hyperref}
\usepackage[capitalize]{cleveref}
\usepackage{indentfirst}

\allowdisplaybreaks
\RequirePackage[colorlinks,citecolor=blue,urlcolor=blue]{hyperref}

\usepackage{subcaption} 
\usepackage{caption}

\newlist{thmlist}{enumerate}{1}
\setlist[thmlist]{label=\alph{thmlisti}., ref=\thetheorem.\alph{thmlisti},noitemsep}

\newlist{deflist}{enumerate}{1}
\setlist[deflist]{label=\alph{deflisti}., ref=\thedefinition.\alph{deflisti}.,noitemsep}

\numberwithin{equation}{section}

\startlocaldefs
\newcommand{\pb}[0]{\mathrm{P}}
\newcommand{\ex}[0]{\mathrm{E}}
\newcommand{\ind}[0]{\mathbf{1}}

\DeclareMathOperator*{\argmax}{arg\,max}

\newlength\myindent
\newcommand\numberthis{\addtocounter{equation}{1}\tag{\theequation}}
\theoremstyle{plain}

\newtheorem{theorem}{Theorem}[section]
\newtheorem{lemma}[theorem]{Lemma}
\newtheorem{corollary}[theorem]{Corollary}
\theoremstyle{definition}

\Crefname{theorem}{Theorem}{Theorems}

\endlocaldefs

\begin{document}

\begin{frontmatter}
\title{Semiparametric Bernstein-von Mises Phenomenon via Isotonized Posterior in Wicksell's problem}

\runtitle{Semiparametric BvM via Isotonized Posterior in Wicksell's problem}

\hypersetup{
pdftitle={Semiparametric Bernstein-von Mises Phenomenon via Isotonized Posterior in Wicksell's problem - F. Gili, G. Jongbloed, A. van der Vaart},
pdfsubject={},
pdfauthor={Francesco ~Gili, Geurt ~Jongbloed, Aad ~van der Vaart},
pdfkeywords={Wicksell's problem, Bayesian Isotonic estimation, Projection Posterior, Dirichlet Process, Bernstein-von Mises},
}

\begin{aug}
\author[A]{\fnms{Francesco}~\snm{Gili} \ead[label=e1]{F.Gili@tudelft.nl}},
\author[A]{\fnms{Geurt}~\snm{Jongbloed}\ead[label=e2]{G.Jongbloed@tudelft.nl}}
\and
\author[A]{\fnms{Aad}~\snm{van der Vaart}\ead[label=e3]{A.W.vanderVaart@tudelft.nl}}
\address[A]{Delft Institute of Applied Mathematics (DIAM), Delft University of Technology, \\ \printead[presep={ \ }]{e1,e2,e3}}

\end{aug}

\begin{abstract}
In this paper, we propose a novel Bayesian approach for nonparametric estimation in Wicksell's problem. This has important applications in astronomy for estimating the distribution of the positions of the stars in a galaxy given projected stellar positions and in materials science to determine the 3D microstructure of a material, using its 2D cross sections. We deviate from the classical Bayesian nonparametric approach, which would place a Dirichlet Process (DP) prior on the distribution function of the unobservables, by directly placing a DP prior on the distribution function of the observables. Our method offers computational simplicity due to the conjugacy of the posterior and allows for asymptotically efficient estimation by projecting the posterior onto the \( \mathbb{L}_2 \) subspace of increasing, right-continuous functions. Indeed, the resulting Isotonized Inverse Posterior (IIP) satisfies a Bernstein-von Mises (BvM) phenomenon with minimax asymptotic variance \( g_0(x)/2\gamma \), where \( \gamma > 1/2 \) reflects the degree of Hölder continuity of the true cdf at \( x \). Since the IIP gives automatic uncertainty quantification, it eliminates the need to estimate \( \gamma \). Our results provide the first semiparametric Bernstein–von Mises theorem for projection-based posteriors with a DP prior in inverse problems.
\end{abstract}



\end{frontmatter}

\section{Introduction}
In the field of stereology, scientists study three-dimensional properties of materials and objects by interpreting their two-dimensional cross-sections. This approach allows the estimation of certain quantities without making costly 3D reconstructions. Consider now these two stereological problems. The first one: astronomers are interested in the distribution of stars in the universe. Globular clusters --- spherical aggregations of stars held together by gravity (cf.\ \cite{33,19}) --- are a topic of particular focus. But how to determine this distribution when given only a 2D photo of the clusters? The second one: having a material that presents a globular microstructure that cannot be scanned (or only at high cost, e.g.\ steel) but only sectioned, how can we determine the distribution of particle sizes inside the material when given only a limited number of 2D cross sections (cf. \cite{49,34})? Despite their seemingly distant nature, these problems share a common thread: the inherent structure of an inverse problem, originally conceptualized by Wicksell (cf.\ \cite{27,1,15,14}).

\begin{figure}[ht]
  \centering
  \includegraphics[width=4.1cm]{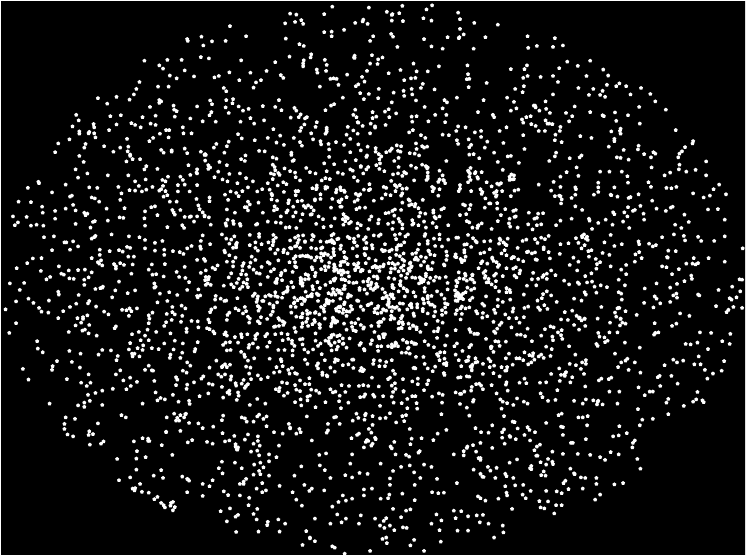}
  \includegraphics[width=5.1cm]{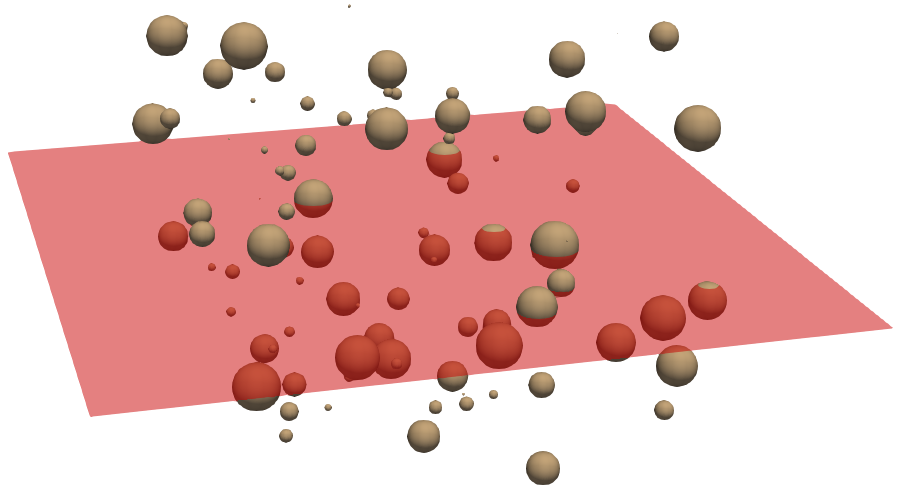}
  \includegraphics[width=2.7cm]{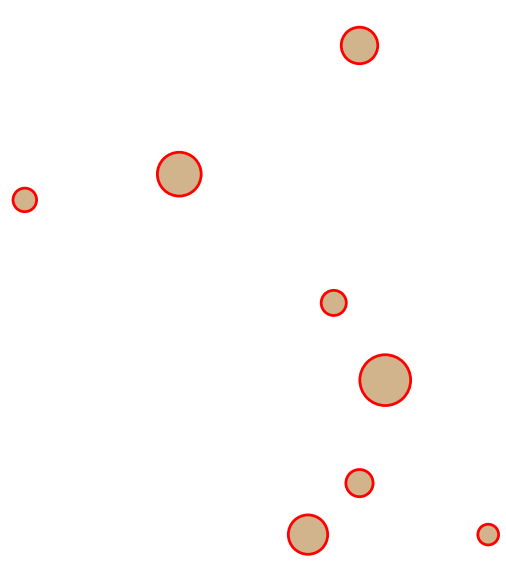}
  \setlength{\belowcaptionskip}{-10pt}
  \caption{Two examples of sampling in stereological inverse problems. Left panel: The Messier 4 (also designated NGC 6121) is a globular cluster in the constellation of Scorpius (c.f.\ \cite{33}). Center and right panel: simulation 3D configuration with spheres and the intersecting plane in material microstructures (from \cite{27}).}\label{fig: 1 sampling}
\end{figure}

The setting of the opaque material having globular microstructure is formalized by assuming that the embedded particle squared radii are i.i.d.\ realizations from a random variable $X^* \sim F^*_0$ and our observations consist of cross sections of the medium (i.e.\ "planes cutting the medium"). Note that the spheres with bigger radii are more likely to be sampled. Thus the distribution of the squared radii of the spheres actually cut is not $F^*_0$ but a biased version of it; we write for the squared radius of the sphere actually cut: $X \sim F_0$. Since we assume that the samples are gathered uniformly, by the Pythagorean Theorem, the observable squared circle radius $Z$ is related to $X$ via: $ Z = (1-U^2) X = Y X,$ where $U \sim U[0,1]$ so $Y \sim \operatorname{Be}(1,1/2)$, with density $y \mapsto \frac{1}{2} (1-y)^{-1/2} \ind_{(0,1)}(y)$. So, the observed squared radii $Z$ are a random fraction of the true squared radii depending on where the plane cuts through the sphere. Consequently, since both $X$ and $Y$ are nonnegative, a version of the density $g_0$ of $Z$ can, for $z\geq 0$, be written as:
\begin{align}\label{eq: g_F_0 expression}
	g_0 (z) = \frac{d}{dz} \int_{(0,\infty)} \pb(Y \leq z/x ) \, dF_0(x) =\int_{z}^{\infty} \frac{dF_0(x)}{2\sqrt{x^2 - xz}}
\end{align}
The problem in the astronomical setting can be modelled as follows: let $\mathbf{Z}= (Z_{p_1}, Z_{p_2},Z_{p_3})$ be a spherically symmetric random vector (the three-dimensional position of a star in a galaxy) of which only $(Z_{p_1}, Z_{p_2})$ can be observed. Under the assumption that the joint of $\mathbf{Z}$ is a function of the distance to the center of the galaxy, astronomers are interested in estimating the distribution of the squared distance of a star to the center of the globular cluster or galaxy $X = Z^2_{p_1}+ Z^2_{p_2}+Z^2_{p_3} \sim F_0$ using $Z = Z^2_{p_1}+ Z^2_{p_2}$. From the assumed spherical symmetry, $\mathbf{Z}
\sim O \mathbf{Z}$ for all orthogonal matrices $O \in \mathrm{O}(3)$. Therefore $\mathbf{U} := \mathbf{Z} / \|\mathbf{Z}\|$ (note $X = \|\mathbf{Z}\|^2$) is uniform on the unit sphere $\mathbb{S}^2$. It follows that the random variable $Z_{p_3}$ (the third coordinate) is distributed as: $\sqrt{X} \cdot U_3$ where $U_3 \sim U[-1,1]$ independent of $X$. Note also that $Y \equiv 1-U_3^2 \sim \operatorname{Be}(1,1/2)$ as above. Because:
\begin{align*}
	\pb (Z \leq z) = \pb ( Z^2_{p_1}+ Z^2_{p_2}+Z^2_{p_3} \leq z + Z^2_{p_3}) = \pb (X \leq z + X \cdot U_3^2) = \pb (YX \leq z),
\end{align*}
it follows that both for the application in materials science and the one in astronomy, the density of the observations $Z \sim G_0$ is the one given on the right-hand side of \eqref{eq: g_F_0 expression}.

By using Abel-integral equation's theory, we can invert \eqref{eq: g_F_0 expression} and obtain (c.f. \cite{19}):
\begin{align}\label{eq: solution range equation}
  &F_0(x) = 1 -\frac{2}{\pi} \sqrt{x} \, V_0(x) - \frac{2}{\pi} \int_{x}^{\infty} \frac{V_0(s)}{2 \sqrt{s}} \, ds,  
\end{align}
where:
\begin{align}\label{eq: V_0}
    V_0(x) \vcentcolon=\int_x^{\infty} \frac{g_0(z)}{\sqrt{z-x}} \: d z.
\end{align}
(The proof of \eqref{eq: solution range equation} is provided in Lemma 1.1 in the Online Supplement \cite{64}). In Corollary \ref{thm: BvM for F}, we formally show that the functional $x \mapsto - \frac{2}{\pi} \int_{x}^{\infty} V_0(s) / (2 \sqrt{s}) \, ds$ can be estimated at the rate $\sqrt{n}$, whereas the functional $V_0(x)$ can be estimated (nonparametrically) only at the slower rate $\sqrt{n}/\sqrt{\log{n}}$. This can be intuitively seen because, by plugging \eqref{eq: V_0} into the second term of \eqref{eq: solution range equation} and using Fubini's theorem to exchange the order of integration:
\begin{align}\label{eq: second equation range}
\int_{x}^{\infty} \frac{V_0(s)}{2 \sqrt{s}} \, ds =\frac{\pi}{2}-\int_0^{\infty} \sin ^{-1} \sqrt{1 \wedge \frac{x}{z}} \,  g_0(z)\, dz 
\end{align}
and $\mathrm{Var}[\sin ^{-1} \sqrt{1 \wedge\left(x / Z\right)}] < \infty$, for $Z \sim G_0$, in contrast with $\mathrm{Var}[(Z-x)^{-1/2}] = \infty$. It thus follows that the estimation of $V_0$ and $F_0$ are essentially equivalent. Furthermore, as shown in \cite{1}, the relation
$ 
F_0^*(x) = 1 - \frac{V_0(x)}{V_0(0)}
$ 
holds, which implies that also in the materials science setting, the primary target of estimation is $V_0$. Accordingly, our focus in this paper will first be on estimating $V_0(x)$, and we will subsequently show that these results naturally extend to the estimation of $F_0(x)$. The latter will be the main object of interest in the astronomical application discussed in the simulation study.

We aim to develop an easy-to-compute and statistically efficient nonparametric Bayesian method for estimating \( V_0(x) \) at a given point \( x \geq 0 \). We begin by placing a Dirichlet Process prior directly on the distribution function of the observables: \( G \sim \operatorname{DP}(\alpha) \). Given an observed sample \( Z_1, \dots, Z_n \) from cdf \( G_0 \) (associated with density $g_0$), the posterior distribution is (c.f.\ \cite{59}, Theorem 4.6 in \cite{4}):
\begin{align}\label{eq: dirichlet posterior}
	G \mid Z_1,\ldots,Z_n \sim \text{DP}(\alpha + n \mathbb{G}_n), \quad \quad  \mathbb{G}_n := \frac{1}{n} \sum_{i=1}^n \delta_{Z_i}.
\end{align}
This choice offers computational simplicity due to the conjugacy of the posterior. Leveraging the inverse relations inherent to the problem, we can estimate the functional of interest using \eqref{eq: solution range equation} by
\begin{align}\label{eq: naive bayes}
	V_{\scriptscriptstyle{G}}(x) := \int_x^\infty \frac{dG(z)}{\sqrt{z - x}}.
\end{align}
This initial step yields the \textit{Naive Bayes Posterior} (NBP). We show below that this satisfies a BvM phenomenon with asymptotic variance \( g_0(x) \). However, the realisations \( V_{\scriptscriptstyle{G}} \) of this plug-in posterior are not necessarily monotone functions. To remedy this, we project them into the \( \mathbb{L}_2 \) subspace of decreasing and right-continuous functions. It turns out that this not only corrects the non-monotonicity, but also improves the estimator. We show that the resulting \textit{Isotonized Inverse Posterior} (IIP) $\hat{V}_{\scriptscriptstyle{G}}(x)$ satisfies a BvM theorem with improved variance that depends on the local smoothness of $V_0$. The smoothness at $x$ is measured through a parameter $\gamma > \frac{1}{2}$ that can be loosely understood as the degree of Hölder continuity at $x$. Using \eqref{eq: solution range equation}, we can define the IIP to target the estimation of $F_0(x)$, which we will denote by $\hat{F}_{\scriptscriptstyle{G}}(x)$. The formal definitions of $\hat{V}_{\scriptscriptstyle{G}}(x)$ and $\hat{F}_{\scriptscriptstyle{G}}(x)$, being quite technical, are given in \eqref{eq: IIP} and \eqref{eq: IIP F0} below, respectively. 

Thus the approach in this paper has three advantages:
\begin{itemize}
	\item[(i)] \textbf{Computational speed}: The Dirichlet posterior \eqref{eq: dirichlet posterior} is explicit and its projection can be effectively computed by the Pool-Adjacent-Violators Algorithm (PAVA, c.f.\ \cite{61}).  
	\item[(ii)]  \textbf{Efficiency and adaptation}: The IIP achieves the asymptotic variance of the efficient estimators introduced in \cite{27} and \cite{10} and it is adaptive to the level of local smoothness of $V_0$.
	\item[(iii)] \textbf{Uncertainty quantification}: The Bayesian framework inherently provides uncertainty quantification without the need to estimate \( \gamma \) or other parameters. We show below that this is asymptotically correct in the frequentist sense.
\end{itemize}

The Naive Bayes Posterior (NBP) is analogous to the naive estimator introduced in \cite{1} and \cite{15} which is constructed as follows: given a sample $Z_1, \ldots, Z_n$ from the density $g_0$ and for $\mathbb{G}_n$ the empirical distribution function based on this sample,
\begin{align}\label{eq: naive estimator of V}
    V_n(x) \vcentcolon=\int_x^{\infty} \frac{d\mathbb{G}_n(z)}{\sqrt{z-x}}.
\end{align}
In \cite{27}, the authors showed that the projection of the naive estimator onto the set of monotone functions, the \textit{Isotonized Inverse Estimator} (IIE) $\hat{V}_n$, is asymptotically efficient and the attained asymptotic variance depends on the local smoothness of $V_0$ in the same way as for the IIP.

Projecting the naive posterior distribution departs from a classical Bayesian framework which might assign a Dirichlet Process (DP) prior to \( F \) and employ a hierarchical model to sample from the density of the observations. A Bernstein-von Mises (BvM) result for this approach is currently unavailable. Moreover, the available algorithms for posterior sampling rely on Gibbs sampling procedures, which are much slower than the method we propose here. 

From a practical perspective, the Bernstein–von Mises (BvM) result established in this paper assures practitioners that the Bayesian methodology based on DP prior provides all the efficiency benefits of the Isotonic Inverse Estimator (IIE) of \cite{27}, but without requiring the estimation of difficult-to-access nuisance quantities such as the density of the observations $g_0(x)$ or the local Hölder smoothness $\gamma > 1/2$. In applications, where one aims to recover the distribution of squared distances of stars from their projected two-dimensional positions or the distribution of the embedded particles' squared radii, this is especially valuable: the underlying smoothness of these distributions (and thus of $V_0$) may vary across radii, and the BvM result holds for all local smoothness levels $\gamma > 1/2$. This significantly improves over existing bootstrap-based uncertainty guarantees of \cite{19}, which only hold under the stronger assumption that $V_0$ is differentiable at the estimation point $x$ (which is even stronger than $\gamma = 1$). 

\subsubsection*{Motivation and connections with the literature} 

There is a vast literature on Wicksell's problem, see for instance \cite{18, 20, 27, 17, 10,1,  2,16,3,19, 9, 15}. Wicksell's problem has present-day applications, where it is common to see the use of the so-called Saltykov method \cite{49}, and numerical discretization methods \cite{34}, often lacking theoretical guarantees. These methods are far from being efficient and in some cases not even consistent. The estimation of $V_0$ given data from $g_0$ is of interest not only for stereological procedures, but also from the mathematical point of view, both because of the unusual rates of convergence and because of the non-standard efficiency theory.

In the nonparametric frequentist approach to Wicksell's problem, the standard rate of convergence for pointwise asymptotic normality in estimating \( V_0 \) is \( (\log n)^{-1/2} \sqrt{n} \) (cf.\ \cite{27, 1, 15}). The first time the naive plug-in estimator was studied is in \cite{15}. Later in \cite{1} the proof of the asymptotic normality of this estimator was simplified considerably. However, viewed as a function, the naive plug-in estimator is severely ill-behaved due to infinite discontinuities and non-monotonicity, and this led to the definition of the IIE in \cite{27,1}. 
The Nonparametric Maximum Likelihood Estimator is studied in \cite{18,3}. For Wicksell's problem, the usual definition of the NPMLE is also ill-posed, as the log-likelihood can attain value infinity. However, work has been done, in \cite{18} and \cite{3}, to obtain likelihood based estimators. Confidence interval results can be found in \cite{20,19}. The results of \cite{20} provide confidence bands by developing a methodology for isotonization procedures, which often lead to nuisance parameters in the limiting distributions, that tackles this issue. However, they do not cover the setting considered in the present paper. Extending their results to our case appears non-trivial, as their proof techniques rely on at least the differentiability of the estimated object—an assumption we do not make. The bootstrap results in \cite{19} are more directly relevant and will be discussed in detail in the simulation section, including a comprehensive numerical comparison. The main advantage of the methodology proposed in this paper, if compared to the bootstrap methods of \cite{19}, is that it relies on significantly less stringent assumptions. In particular, our results do not depend on the assumption that the underlying distribution function $F_0$ is differentiable at $x$. Instead, our methodology is provably consistent for all local smoothness levels $\gamma > 1/2$ of $V_0$, where $\gamma$ is defined in condition \eqref{eq: condition roughness} below.

In \cite{10}, the authors investigate efficient estimation in the classical Wicksell problem, focusing on the scenario where the distribution function $F^*_0$ of interest belongs to a Hölder class of smoothness with degree $\gamma > \frac{1}{2}$. They propose a kernel density-based estimator that, under appropriate bandwidth selection, achieves asymptotic efficiency. However, a notable limitation is that the estimator requires prior knowledge or guesstimation of the smoothness parameter $\gamma$, making it non-adaptive. Moreover, efficiency is only guaranteed under global smoothness assumptions, i.e., the entire distribution function must belong to a global Hölder class. Despite these limitations, their work remains the first foundational contribution to semiparametric efficiency in Wicksell's problem. \cite{27} builds upon earlier work of \cite{10} and \cite{1} by addressing both the adaptivity and locality of efficient estimation in Wicksell's problem. The key innovation consists in showing that the Isotonized Inverse Estimator (IIE) achieves asymptotic efficiency under the same local smoothness assumptions on the target function $F^*_0$ of the present paper. Specifically, instead of assuming that $F^*_0 \in \Lambda_\gamma$ globally as in \cite{10}, the analysis in \cite{27} only requires that $F^*_0$ satisfies local Hölder-type smoothness at the evaluation point $x$ and at the origin. The estimator adapts automatically to the local regularity of $F^*_0$, offering a more flexible framework for practitioners. Additionally, \cite{27} provides a local asymptotic minimax lower bound, showing that the IIE attains the asymptotically optimal risk in this refined setting. The estimator thus achieves the same efficiency as the KDE-based approaches, but without requiring tuning of smoothness parameters, and without imposing strong global assumptions. Finally, other results concerning efficient estimation in Wicksell's problem are contained in \cite{50} and \cite{17}. \cite{50} shows that in the case the underlying cdf is constant on a \textit{known} interval of strictly positive length the IIE evaluated at $x$ is not efficient, for $x$ in the interval of constancy. However, it behaves in practice very closely to informed estimators that make use of the fact that the cdf is constant on that interval. Finally, one of the authors of \cite{10} did some further work in \cite{17} where it was shown how to attain efficient smoothing for the plug-in estimator under squared losses. 

While Bayesian methods for inverse problems are widely adopted, their theoretical guarantees remain scarce, particularly in infinite-dimensional settings. Although recent years have seen significant advances in this area \cite{62,35,37,41,36,44,38,39,40,42,26,43}, the technical frameworks developed there do not readily apply to the specific inverse problem and Dirichlet process prior setting of this paper. We provide a detailed discussion of this related work in Supplement \cite{64}. A similar methodology to the one proposed in this paper --- using an unconstrained conjugate prior for the distribution function of the observables and projecting the posterior samples into the space of monotone functions --- has been very recently adopted in isotonic regression settings (c.f.\ \cite{54,55,56}), isotonic density estimation (c.f.\ \cite{60}) and ODEs settings (c.f.\ \cite{57,58}). \cite{57,58} develop Bayesian two-step methods for parameter inference in ODE models by embedding them in a nonparametric regression framework with B-spline priors. The parameter posterior is obtained by minimizing a discrepancy between the estimated regression function and the ODE solution, either via derivative matching or Runge–Kutta approximations. \cite{58} further establishes a Bernstein–von Mises theorem for the Runge–Kutta-based method, showing that the resulting Bayes estimator is asymptotically efficient. However, their framework relies on a very different type of prior modelling and the projection method does not entail isotonization. In \cite{54}, the authors study a nonparametric Bayesian regression model for a response variable \( Y \) with respect to a predictor variable \( X \in [0,1] \), given by \( Y = f(X) + \varepsilon \), where \( f \) is a monotone increasing function and \( \varepsilon \) is a mean-zero random error. They propose a "projection-posterior" estimator \( f^* \), constructed using a finite random series of step functions with normal basis coefficients as a prior for \( f \). Their main finding is that when \( f^* \) is centered at the maximum likelihood estimator (MLE) and rescaled appropriately, the Bernstein-von Mises (BvM) theorem does not hold. However, when \( f^* \) is instead centered at the true function \( f_0 \), a BvM result is obtained. This phenomenon resonates with the known inconsistency of the bootstrap for the Grenander estimator in isotonic regression (c.f.\ \cite{63}) and represents a key difference with the present work. In our case, the "projection-posterior" IIP \( \hat{V}_{\scriptscriptstyle{G}} \), when centered at the IIE \( \hat{V}_n \), is shown to converge in distribution.  This may also not be surprising given the existing results in \cite{19}, where bootstrap methods are shown to be consistent for Wicksell’s problem. 

One of the key advantages our approach shares with the works in \cite{57,58,54,55,56} is the computational efficiency, driven by posterior conjugacy, as well as the relative simplicity of the asymptotic analysis. However, unlike in \cite{54}, our method ensures that the limiting distribution in the Bernstein-von Mises (BvM) theorem coincides with that of the asymptotically efficient estimator --- the IIE of \cite{27}.  

\section{Main results}\label{sec: main results}

We start by proving a BvM result for the \textit{Naive Bayes Posterior} (NBP), analogous to the result obtained for the naive estimator (c.f.\ Theorem 2 in \cite{1}). For our first result, we need our version of $g_0$ to be well-defined and bounded, i.e.:
\begin{align*}
    \| g_0 \|_{\infty} := \sup_{z \geq 0} \bigg| \int_{z}^{\infty} \frac{dV_0(s)}{\pi \sqrt{s - z}} \bigg| < \infty. \numberthis \label{eq: condition finite g at x}
\end{align*}
Assumption \eqref{eq: condition finite g at x} asserts that projected 2D star positions do not become arbitrarily concentrated. From a mathematical perspective, without this assumption the inverse problem becomes highly unstable: small errors in the data can cause arbitrarily large errors in estimating $V_0$. For instance, if $g_0$ is unbounded at $x$, this leads to instability in estimating $V_0(x)$, reflected in an infinite asymptotic variance. Thus, boundedness of $g_0$ ensures stability of the inverse mapping. Mathematically, if $F_0$ is Hölder continuous with degree $\gamma > 1/2$, then $g_0$ is bounded (cf.\ Lemma A.3 in \cite{27}).

In equation \eqref{eq: claim naive} and \eqref{eq: claim} below we use the notion of conditional weak convergence in distribution (in probability). For its precise definition we refer the reader to \eqref{eq: formal definition conditional convergence}-\eqref{eq: formal definition conditional convergence 1} in Appendix A.
\begin{theorem}[BvM for NBP]\label{thm: BvM naive}
	Let $x \geq 0$ and suppose that \eqref{eq: condition finite g at x} holds true. Let $G \sim \operatorname{DP}(\alpha)$, with $\int_{x}^{\infty} (s-x)^{-1/2} \, d \alpha(s) < \infty$. Then, in probability under $G_0$: 
	\begin{align}\label{eq: claim naive}
	\frac{\sqrt{n}}{\sqrt{\log n}}\big( V_{\scriptscriptstyle{G}}(x)-V_n(x)\big) \mid Z_1, \ldots, Z_n \rightsquigarrow N\big(0, g_0(x)\big).
	\end{align}
\end{theorem}
\begin{proof}
	For a measure $P$ and a function $f$, we denote by $Pf$ the integral $\int f \, dP$, this holding also for empirical measures like $\mathbb{G}_n = \frac{1}{n} \sum_{i=1}^n \delta_{Z_i}$. 
	
	Since $G \mid Z_1,\ldots,Z_n \sim \text{DP}(\alpha + n \mathbb{G}_n)$ by proposition G.10 in \cite{4}, for \( Q \sim \text{DP}(\alpha) \), \( \mathbb{B}_n \sim \text{DP}(n \mathbb{G}_n) \), \( V_n \sim \text{Be}(|\alpha|, n) \) and $f(z) = 1/\sqrt{(z-x)_+} $, $x \geq 0 $, conditionally on $Z_1,\ldots,Z_n$, we have the following representation in distribution:
	\begin{align}\label{eq: first representation in distribution}
	G f = V_n Q f + (1 - V_n) \mathbb{B}_n f,
	\end{align}
	where \( Q  \), \( \mathbb{B}_n \), \( V_n  \) are independent. First note that by assumption $Qf$ is a.s.\ well-defined (c.f.\ Remark 4.4 in \cite{4}). By Proposition G.2 in \cite{4} for \( \varepsilon_i \overset{\text{i.i.d.}}{\sim} \text{Exp}(1) \) and \( W_{ni} := \varepsilon_i / \bar{\varepsilon}_n \), we have \( (W_{n1}, \ldots, W_{nn}) \sim \text{Dir}(n; 1/n, \dots, 1/n) \). Since \( \mathbb{B}_n \sim \text{DP}(n \mathbb{G}_n) \), for $\delta_{Z_i}$ Dirac measure at $Z_i$, we have the following representation in distribution of $\mathbb{B}_n$ (c.f.\ Example 3.7.9 in \cite{5}):
	\[
	\mathbb{B}_n = \frac{1}{n} \sum_{i=1}^n W_{ni} \delta_{Z_i}.
	\]
	Because of this representation (c.f.\ section 3.7.2. in \cite{5}), $\mathbb{B}_n$ is also known as \textit{Bayesian bootstrap}. Conclude that: 
	\[
		\frac{\sqrt{n}}{\sqrt{\log n}} (\mathbb{B}_n - \mathbb{G}_n) f = \frac{1}{\sqrt{n \log n}} \sum_{i=1}^n \left( W_{ni} - 1 \right) f(Z_i).
	\]
	 Because $|\bar{\varepsilon}_n - 1| = O_p(1/\sqrt{n})$, we study: $(n \log n)^{-1/2} \sum_{i=1}^n \left( \varepsilon_i - 1 \right) f(Z_i)$. Now note that by Lemma 2.4 in the Online Supplement \cite{64}:
	\begin{align}\label{eq: weak conv for NBP}
		\frac{1}{\sqrt{n \log{n}}} \sum_{i=1}^n \frac{(\varepsilon_i -1)}{\sqrt{Z_i - x}} \ind_{\{ Z_i >x \}} \mid Z_1, \ldots, Z_n \rightsquigarrow N\left(0, g_0(x)\right) \, \, \text { under } \, G_0 \text {. }
	\end{align}
	Because $Qf$ is a.s.\ well-defined, we obtain the claim by Slutsky's lemma upon noticing $ \frac{\sqrt{n}}{\sqrt{\log{n}}} V_n Q f \stackrel{G_0}{\rightarrow} 0$.
\end{proof}
We consider a projection-type posterior by isotonizing draws from the obtained naive posterior \eqref{eq: naive bayes} as follows. First define, for $G \sim \text{DP}(\alpha + n \mathbb{G}_n)$ and $\mathbb{G}_n$, $G_0$ as above:
\begin{align*}
&U_{\scriptscriptstyle{G}}(x) := \int_{0}^x V_{\scriptscriptstyle{G}}(y) \, dy = 2 \int_{0}^\infty \sqrt{z} \, dG(z) - 2 \int_x^\infty \sqrt{z - x} \, dG(z), \numberthis \label{eq: U_G}\\
&U_n(x) := \int_{0}^x V_n(y) \, dy = 2 \int_{0}^\infty \sqrt{z} \, d\mathbb{G}_n(z) - 2 \int_x^\infty \sqrt{z - x} \, d\mathbb{G}_n(z),\numberthis \label{eq: U_n} \\
&U_0(x) := \int_{0}^x V_0(y) \, dy = 2 \int_{0}^\infty \sqrt{z} \, dG_0(z) - 2 \int_x^\infty \sqrt{z - x} \, dG_0(z). \numberthis \label{eq: U_0}
\end{align*}

Furthermore, for a continuous function $k : [0,\infty) \to [0,\infty)$, the Least Concave Majorant (LCM) and its right-hand side derivative, are given for $x \in [0,\infty)$ by:
\begin{align*}
    &k^*(x) := \inf {\{ f(x) \, : \, f \geq k \, \, \text{on} \, \, [0,\infty), \, f \, \text{concave} \, \, \text{on} \, \, [0,\infty) \}},\\
    & (k^*)^{\prime}_+ (x):= \lim_{h \to 0^+} \frac{k^*(x + h) - k^*(x)}{h}.
\end{align*}
We denote by $\hat{V}_n$ (IIE) the right-hand side derivative of $U^*_n$, the LCM of $U_n$. Similarly, $\hat{V}_{\scriptscriptstyle{G}}$ is the right-hand side derivative of $U^*_{\scriptscriptstyle{G}}$, the LCM of $U_{\scriptscriptstyle{G}}$, i.e., the isotonization of a draw from the naive Bayes posterior. Thus, for $G$ as in \eqref{eq: dirichlet posterior}:  
\begin{align*}
	&\hat{V}_n(x) := (U^*_n)^{\prime}_+ (x), \numberthis \label{eq: IIE}\\
	&\hat{V}_{\scriptscriptstyle{G}}(x) :=  (U^*_{\scriptscriptstyle{G}})^{\prime}_+(x). \numberthis \label{eq: IIP}
\end{align*} 
By Lemma 2 in \cite{1}, $\hat{V}_n$ and $\hat{V}_{\scriptscriptstyle{G}}$ can be viewed as (equivalent to) $\mathbb{L}_2$-projections of the naive estimators $V_n$ and $V_{\scriptscriptstyle{G}}$ into the space of nonnegative, decreasing and right-continuous functions on $[0,\infty)$. We can finally define the Naive posterior and the IIP to target the estimation of $F_0$ as follows:
\begin{align*}
    &F_{\scriptscriptstyle{G}}(x) := 1 -\frac{2}{\pi} \sqrt{x} \, V_{\scriptscriptstyle{G}}(x) - \frac{2}{\pi} \int_{x}^{\infty} \frac{V_{\scriptscriptstyle{G}}(s)}{2 \sqrt{s}} \, ds, \numberthis \label{eq: naive post F0} \\
	&\hat{F}_{\scriptscriptstyle{G}}(x) := 1 -\frac{2}{\pi} \sqrt{x} \, \hat{V}_{\scriptscriptstyle{G}}(x) - \frac{2}{\pi} \int_{x}^{\infty} \frac{\hat{V}_{\scriptscriptstyle{G}}(s)}{2 \sqrt{s}} \, ds,  \numberthis \label{eq: IIP F0}
\end{align*}
and analogously we define $F_n(x)$, with $V_n$ instead of $V_{\scriptscriptstyle{G}}$ in \eqref{eq: naive post F0} and $\hat{F}_n(x)$, with $\hat{V}_n$ instead of $\hat{V}_{\scriptscriptstyle{G}}$ in \eqref{eq: IIP F0}.

We now state two additional conditions needed for Theorem \ref{thm: BvM}. First, a natural assumption is a moment constraint for $Z$, (which can be translated into a moment condition on $X$ c.f.\ \cite{27}): 
\begin{align*}
    \quad \int_{0}^{\infty} z \, g_0(z) \, dz < \infty. \numberthis \label{eq: finite first moment}
\end{align*}
This is a mild assumption, which is practically always satisfied, given that materials studied in materials science are of bounded size and globular clusters typically do not extend indefinitely. Additionally, assume $\exists$ $K >0$ and $\gamma > \frac{1}{2}$ such that
\begin{align*}
    &  \int_{0}^{1} ( V_0(x) - V_0(x+u \delta) ) \, du \sim \operatorname{sgn}(\delta)|\delta|^{\gamma} K, \quad \quad \text{as} \: \delta \rightarrow 0, \numberthis \label{eq: condition roughness} 
\end{align*}
where $\sim$ means that the ratio of the left-hand side and of the right-hand side goes to 1 as $\delta \rightarrow 0$. The constant $\gamma$ indicates the level of smoothness at $x$ of $V_0$. This is analogous to equation (5) from \cite{27}. Note that it is possible to see what this condition means in terms of the underlying $F_0$, using \eqref{eq: solution range equation}, we obtain that $\exists$ $K >0$ and $\gamma > \frac{1}{2}$ such that as $\delta \rightarrow 0$:
\begin{align}\label{eq: condition roughness F0}
	\int_{0}^{1} ( \sqrt{x} \, V_0(x) - \sqrt{x+ u \delta} \, V_0(x+u \delta) ) \, du \sim  \operatorname{sgn}(\delta)|\delta|^{\gamma} K - \frac{V_0(x)}{4 \sqrt{x}} \operatorname{sgn}(\delta)|\delta|, 
\end{align}
(upon using a Taylor expansion for $\int_0^1 (\sqrt{x} - \sqrt{x+ u \delta}) \, du$). Because $\frac{2}{\pi} \int_{x}^{\infty} \frac{V_0(s)}{2 \sqrt{s}} \, ds$ is Lipschitz continuous, the local smoothness of $F_0$ at $x$ is $\min(\gamma, 1)$. Therefore for the most relevant levels of smoothness, i.e.\ $\gamma \in (\frac{1}{2},1]$, the smoothness of $V_0$ carries over to $F_0$. Because $ 
F_0^*(x) = 1 - \frac{V_0(x)}{V_0(0)}
$, equation (5) from \cite{27} is equivalent to \eqref{eq: condition roughness} (note that in that paper $F$ denotes the current $F_0^*$). Thus, for a better mathematical understanding of the role of $\gamma$, we refer the interested reader to the discussion given in \cite{27}, in particular Example 2.1, 2.2, 2.3 and Remark 2.1. 

From a physical point of view, the role of $\gamma$ is also interpretable. In a globular cluster, different regions (e.g., dense core vs. sparse halo) may correspond to different levels of smoothness of $F_0$, and the semiparametric theory reflects this. The local adaptivity of the Bayesian estimator — encoded in the limiting variance of the Bernstein–von Mises — mirrors this phenomenon. If there is a sharp increase in the concentration of stars with size close to $x$ (i.e., $\gamma$ close to $1/2$), then the 2D observations are much noisier. This makes the estimation more difficult and results in a higher asymptotic variance. 

\begin{theorem}[BvM for IIP]\label{thm: BvM}
	Let $x \geq 0$, suppose that \eqref{eq: condition finite g at x}, \eqref{eq: finite first moment}, and \eqref{eq: condition roughness} hold true. Moreover let $G \sim \operatorname{DP}(\alpha)$, where the base measure $\alpha$ has bounded density. Then, in probability under $G_0$: 
	\begin{align}\label{eq: claim}
		\frac{\sqrt{n}}{\sqrt{\log n}}\left(\hat{V}_{\scriptscriptstyle{G}}(x)-\hat{V}_n(x)\right) \mid Z_1, \ldots, Z_n \rightsquigarrow N\left(0, \frac{g_0(x)}{2 \gamma}\right).
	\end{align}
\end{theorem}
\begin{proof}

The proof essentially consists of three parts. The first part is to rewrite the event
$$
\left\{\sqrt{\frac{n}{\log n}}\left(\hat{V}_{\scriptscriptstyle{G}}(x) - \hat{V}_n(x)\right) \leq a\right\}
$$
given $Z_1,\ldots,Z_n$ in terms of the location of the maximum (arg max) of a stochastic process $t \mapsto \tilde{Z}^{\scriptscriptstyle{G}}_n(t)-a t$. It turns out that the process $\tilde{Z}^{\scriptscriptstyle{G}}_n(t)$ can be decomposed into three essential parts. The second step consists in studying these processes. One process can be interpreted as the Bayesian bootstrap counterpart of the process studied in the proof of Theorem 2.1 in \cite{27} (c.f.\ eq.\ (14) in \cite{27}). This process converges conditionally in distribution to a Gaussian process $\mathbb{Z}$, which we show can be written as $\mathbb{Z}(t) = tX$, where $X \sim N(0,g_0(x)/(2\gamma))$. A second process (see \eqref{eq: deterministic part process}), which is a deterministic process, whose behavior is controlled by the local smoothness of $V_0$ at $x$ and which converges to $- |t|^{\gamma+1} K$ (with $\gamma$ and $K$ as in \eqref{eq: condition roughness}). Lastly, a third remainder process which converges uniformly to zero with respect to the topology of uniform convergence on compacta. The final step is to prove that from the conditional convergence in distribution of $t \mapsto \tilde{Z}^{\scriptscriptstyle{G}}_n(t)-a t$ to $\mathbb{Z}(t)-a t$ it follows that the arg max of ($\tilde{Z}^{\scriptscriptstyle{G}}_n(t)-a t$) converges in distribution to the arg max of $\mathbb{Z}(t)-a t$. Since the arg max functional is certainly not continuous on the space of continuous functions equipped with the uniform norm, the simplest form of the continuous mapping theorem cannot be applied. For this reason, we prove in Lemma \ref{lemma: cond argmax}, specifically for this application, a conditional argmax continuous mapping theorem (extension of Theorem 3.2.2 in \cite{5}). From \eqref{eq: conditional stoch boundedness} below, it follows that the sequence of locations of the maximum of ($\tilde{Z}_n(t)-a t$) is conditionally uniformly tight. This enables us to use Lemma \ref{lemma: cond argmax}. Moreover, since the limiting process turns out almost surely to have a unique maximum, the asserted conditional convergence in distribution in \eqref{eq: claim} follows.

For sequences $\delta_n = \sqrt{\log{n}}/\sqrt{n}$ and $\delta_n^{*}=\delta_n^{1 / \gamma}$ where $\gamma >1 / 2$ satisfies condition \eqref{eq: condition roughness}, define
\begin{align}\label{eq: Z_n process}
	Z^{\scriptscriptstyle{G_0}}_n(t) & :=\delta_n^{-1}\left(\delta_n^*\right)^{-1}\left\{U_n\left(x+\delta_n^* t\right)-U_n(x)-U_0\left(x+\delta_n^* t\right)+U_0(x)\right\} \\
	& =2 \delta_n^{-1} (\delta^*_n)^{-1} \int\left\{\sqrt{(z-x)_+} -\sqrt{(z-x-\delta^*_n t)_+} \right\} \, d\left(\mathbb{G}_n-G_0\right)(z).
\end{align}
By Lemma \ref{eq: lemma equiv to claim} in Appendix A it suffices to show, for $X \sim N\left(0, \frac{g_0(x)}{2 \gamma}\right)$ and every $a \in \mathbb{R}$:
\begin{align}\label{eq: equiv to claim}
	\Pi_n\left( \sqrt{\frac{n}{\log n}}\left(\hat{V}_{\scriptscriptstyle{G}}(x)-V_0(x)\right) -Z^{\scriptscriptstyle{G_0}}_n(1) \leq a \mid Z_1, \ldots, Z_n \right) \xrightarrow{G_0} \mathrm{P}(X \leq a).
\end{align}
Define the process \( T_{\scriptscriptstyle{G}} \) for \( a > 0 \):
\[
T_{\scriptscriptstyle{G}}(a) := \inf \{ t \geq 0 : U_{\scriptscriptstyle{G}}(t) - at \text{ is maximal} \}.
\]
By the switch relation (see Lemma 3.2 in \cite{2}):
\[
T_{\scriptscriptstyle{G}}(a) \leq x \iff \hat{V}_{\scriptscriptstyle{G}}(x) \leq a.
\]
Therefore, for sufficiently large \( n \), conditionally on $Z_1, \ldots, Z_n$:
\[
\delta_n^{-1} ( \hat{V}_{\scriptscriptstyle{G}}(x) - V_0(x) ) - Z^{\scriptscriptstyle{G_0}}_n(1) \leq a 
\iff \left( \delta_n^* \right)^{-1} \left( T_{\scriptscriptstyle{G}} \left( V_0(x) + \delta_n (a + Z^{\scriptscriptstyle{G_0}}_n(1)) \right) - x \right) \leq 0.
\]
Because the location of a maximum of a function is equivariant under translations and invariant under multiplication by a positive number and addition of a constant:
\begin{align*}
	&\left( \delta_n^* \right)^{-1} \left( T_{\scriptscriptstyle{G}} \left( V_0(x) + \delta_n (a + Z^{\scriptscriptstyle{G_0}}_n(1)) \right) - x \right)\\
	& = \left( \delta_n^* \right)^{-1} \big( \inf \big\{ x + t \geq 0 : U_{\scriptscriptstyle{G}}(x + t) - V_0(x)(x+t) \\
    &\quad \quad \quad \quad \quad - \delta_n (a + Z^{\scriptscriptstyle{G_0}}_n(1))(x + t) \, \, \text{ is max} \big\} - x \big)\\
	&= \inf \big\{t \geq-(\delta^*_n)^{-1} x: U_{\scriptscriptstyle{G}}\left(x+\delta^*_n t\right)-U_{\scriptscriptstyle{G}}(x) -U_n\left(x+\delta^*_n t\right)+U_n(x)\\
	&\quad + U_0(x+\delta^*_nt) - U_0(x)-V_0(x) \delta^*_n t-\delta_n \delta^*_n (a t + Z^{\scriptscriptstyle{G_0}}_n(1)t - Z^{\scriptscriptstyle{G_0}}_n(t)) \, \text { is max}\big\} \numberthis \label{eq: useful for stoch}\\
	&= \inf \big\{t \geq-(\delta^*_n)^{-1} x: \tilde{Z}^{\scriptscriptstyle{G}}_n(t)-a t \, \text { is max}\big\}.
\end{align*}
where we define the localized process \( \tilde{Z}_n^{\scriptscriptstyle{G}} \) by,
\begin{align*}
\tilde{Z}_n^{\scriptscriptstyle{G}}(t) &:= \delta_n^{-1} \left( \delta_n^* \right)^{-1} \big\{ U_{\scriptscriptstyle{G}} \left( x + \delta_n^* t \right) - U_{\scriptscriptstyle{G}}(x) - U_n \left( x + \delta_n^* t \right) + U_n(x) \big\} \numberthis \label{eq: first components Z tilde}\\
&\quad + \delta_n^{-1} \left( \delta_n^* \right)^{-1} \big\{ U_0 \left( x + \delta_n^* t \right) - U_0(x) - V_0(x) \delta^*_n t \big\} +  Z^{\scriptscriptstyle{G_0}}_n(t) -  t Z^{\scriptscriptstyle{G_0}}_n(1) \numberthis \label{eq: second components Z tilde}
\end{align*}
for $t \in I_x$, where $I_x := [0, \infty) \text{ if } x = 0, \; \text{and} \; I_x := \mathbb{R} \text{ if } x > 0$.

We determine the asymptotic behavior of the process $\tilde{Z}^{\scriptscriptstyle{G}}_n$ by studying separately the components in \eqref{eq: first components Z tilde} and \eqref{eq: second components Z tilde}, and next we conclude by applying Slutsky's lemma. First note that $|Z^{\scriptscriptstyle{G_0}}_n(1)t - Z^{\scriptscriptstyle{G_0}}_n(t)| = o_{\scriptscriptstyle{G_0}}(1)$ uniformly in $t$ in compacta by the proof of Theorem 2.1 in \cite{27} (c.f.\ display after equation (16)). Furthermore, using assumption \eqref{eq: condition roughness} (as in equation (13) in \cite{27}), uniformly in $t$ in compacta:
\begin{align}\label{eq: deterministic part process}
	\delta_n^{-1} \left( \delta_n^* \right)^{-1} \big\{ U_0 \left( x + \delta_n^* t \right) - U_0(x) - V_0(x) \delta^*_n t \big\} =  - |t|^{\gamma+ 1} (K + o(1)).
\end{align}
This determines the behavior of \eqref{eq: second components Z tilde}. Now we study the posterior behaviour of the process in \eqref{eq: first components Z tilde}, which can be interpreted as a \textit{Bayesian bootstrap} counterpart of the process \eqref{eq: Z_n process} studied in the proof of Theorem 2.1 in \cite{27} (c.f.\ eq.\ (14) in \cite{27}):
\begin{align*}
	Z^{\scriptscriptstyle{G}}_n(t) &:= \delta_n^{-1} \left( \delta_n^* \right)^{-1} \big\{ U_{\scriptscriptstyle{G}} \left( x + \delta_n^* t \right) - U_{\scriptscriptstyle{G}}(x) - U_n \left( x + \delta_n^* t \right) + U_n(x) \big\} \\
	&= 2 \delta_n^{-1} \left( \delta_n^* \right)^{-1} \int\left\{\sqrt{(z-x)_+} -\sqrt{(z-x-\delta^*_n t)_+} \right\} \, d\left(G-\mathbb{G}_n\right)(z)
\end{align*}
As in proof of Theorem \ref{thm: BvM naive} in equation \eqref{eq: first representation in distribution}, we have by Proposition G.10 in \cite{4}, for independent \( Q \sim \text{DP}(\alpha) \), \( \mathbb{B}_n \sim \text{DP}(n \mathbb{G}_n) \), and \( V_n \sim \text{Be}(|\alpha|, n) \), the following representation in distribution:
\begin{align}\label{eq: representation process G}
	G f_t^n = V_n Q f_t^n + (1 - V_n) \mathbb{B}_n f_t^n,
\end{align}
where $f_t^n(z) := 2(\delta_n^*)^{-1} \big\{ \sqrt{(z-x)_+} -\sqrt{(z-x-\delta^*_n t)_+}  \big\}$ and thus the following representation in distribution of the process \( Z^{\scriptscriptstyle{G}}_n(t) \):
\[
Z^{\scriptscriptstyle{G}}_n(t) = \delta_n^{-1} V_n (Q - \mathbb{G}_n) f_t^n + \delta_n^{-1} (1 - V_n) (\mathbb{B}_n - \mathbb{G}_n) f_t^n.
\]
This representation is crucial for the asymptotic analysis of the process \( Z^{\scriptscriptstyle{G}}_n(t) \). The first thing to do is to identify which parts of the process \( Z^{\scriptscriptstyle{G}}_n(t) \) do not give any contribution to limiting distribution of the process, as they converge conditionally to zero in probability. The second part is to identify the limiting distribution of the remaining parts of the process, which will be a Gaussian process. In Appendix A, we give the \hyperlink{proof of reduction}{proof} that uniformly in $t$ in compacta:
\begin{align}\label{eq: representation without remainders}
	Z^{\scriptscriptstyle{G}}_n(t) = \mathbb{Z}^*_n (t) + o_p(1),
\end{align}
where for \( \varepsilon_i \overset{\text{i.i.d.}}{\sim} \text{Exp}(1) \), the process $\mathbb{Z}^*_n$ is defined as:
\begin{align}\label{eq: mathbbZ process}
    \mathbb{Z}^*_n (t) := \frac{1}{\sqrt{n \log n}} \sum_{i=1}^n \left( \varepsilon_i - 1 \right) f_t^n(Z_i).
\end{align}
Then by Slutsky's lemma, if we determine the conditional weak limit of \( (\mathbb{Z}_n^\ast(t) \, : \, t \in I_x) \), then we determine the one of \( (Z^{\scriptscriptstyle{G}}_n(t) \, : \, t \in I_x) \). We now study the marginals of $\mathbb{Z}^*_n$. Note that by Lemma 2.5 in the Online Supplement \cite{64} for fixed $t \in I_x$, in probability under $G_0$:
\begin{align}\label{eq: weak conv for IIP}
	\mathbb{Z}^*_n (t) \mid Z_1, \ldots, Z_n \rightsquigarrow N\left(0, t^2  g_0(x)/ (2 \gamma) \right),
\end{align}
and again by Lemma 2.5 in the Online Supplement \cite{64} that:
\begin{align}\label{eq: cov structure}
\text{Cov}\left(\mathbb{Z}_n^\ast(t), \mathbb{Z}_n^\ast(s) \,\middle|\, Z_1, \dots, Z_n\right) 
= \frac{1}{2 \gamma} g_0(x) st \left( 1 - O_{\scriptscriptstyle{G_0}}\left(\frac{\log \log n}{\log n} \right) \right) + O_{\scriptscriptstyle{G_0}}\left(\frac{1}{\log n}\right),
\end{align}
By the conditional Chebyshev inequality and \eqref{eq: cov structure} we conclude $\forall \: \varepsilon >0$:
\begin{align*}
    \pb \left( | \mathbb{Z}^*_n (t) - t \mathbb{Z}^*_n (1) | > \varepsilon \mid Z_1, \ldots, Z_n \right) \leq \frac{1}{\varepsilon^2} \mathrm{Var}( \mathbb{Z}^*_n (t) - t \mathbb{Z}^*_n (1)  \mid Z_1, \ldots, Z_n ) \stackrel{G_0}{\rightarrow} 0.
\end{align*}
By \eqref{eq: representation without remainders} conclude that \eqref{eq: weak conv for IIP} holds true with $\mathbb{Z}^*_n $ replaced by $Z^{\scriptscriptstyle{G}}_n$. Therefore, the finite-dimensional distributions of $Z^{\scriptscriptstyle{G}}_n$ converge conditionally in distribution to the finite-dimensional distributions of the process: $ \mathbb{Z}(t) = tX$, where $X \sim N(0,g_0(x)/(2\gamma))$. 

Lemma \ref{lemma: cond asymp equicontinuity} gives the conditional asymptotic equicontinuity of the process $\mathbb{Z}^*_n$. The conditional asymptotic equicontinuity together with the conditional convergence in distribution of the finite dimensional distributions, implies in probability under $G_0$ as $n \rightarrow \infty$:
\begin{align}\label{eq: weak conv Z^*_n}
	\mathbb{Z}^*_n \mid Z_1,\ldots, Z_n \rightsquigarrow \mathbb{Z} \quad \text{in} \quad \ell^{\infty}(K).
\end{align}
for every compact $K \subset I_x$ and $\ell^{\infty}$ being the space of bounded functions on $K$ equipped with the supremum norm. For a review of the conditional weak convergence in $\ell^{\infty}$, see \cite{5} (c.f.\ Example 1.5.1 and Theorem 1.5.4). By \eqref{eq: representation without remainders} and \eqref{eq: weak conv for IIP} conclude that \eqref{eq: weak conv Z^*_n} holds true with $\mathbb{Z}^*_n $ replaced by $Z^{\scriptscriptstyle{G}}_n$, i.e.\ the conditional weak limit of $Z^{\scriptscriptstyle{G}}_n$ relative to the topology of uniform convergence on compacta is the same Gaussian process. Altogether, recalling also \eqref{eq: deterministic part process}, we can conclude that conditionally on $Z_1,\ldots, Z_n$ the weak limit of $(\tilde{Z}_n^{\scriptscriptstyle{G}}(t) - at \, : \, t \in K )$ in $\ell^{\infty}$ is given by $( tX - |t|^{\gamma+1} K - at \, : \, t \in K )$ for every compact $K\subset I_x$.

In order to apply Lemma \ref{lemma: cond argmax}, notice that \( \mathbb{Z}\) has a.s.\ continuous sample paths and a.s.\ unique maxima (see Lemma 2.6 in \cite{12}) and that $\hat{t}^{\scriptscriptstyle{G}}_n :=  \inf \big\{t \geq-(\delta^*_n)^{-1} x: \tilde{Z}^{\scriptscriptstyle{G}}_n(t)-a t \, \text { is max}\big\}$ is an argmax process conditionally on $Z_1, \ldots, Z_n$, satisfying the second condition on Lemma \ref{lemma: cond argmax}. It is thus enough to show that there exists a conditionally stochastically bounded sequence \((\tilde{M}_n)\) that satisfies:
\begin{align}\label{eq: conditional stoch boundedness}
\Pi_n \left( \left( \delta_n^* \right)^{-1} \left| T_{\scriptscriptstyle{G}} \left( V_0(x) + \delta_n (a + Z^{\scriptscriptstyle{G_0}}_n(1)) \right) - x \right| \leq \tilde{M}_n \mid Z_1, \ldots, Z_n \right) 
\stackrel{G_0}{\to} 1,
\end{align}
The \hyperlink{proof of stoch boundedness}{proof} of \eqref{eq: conditional stoch boundedness} is given in Appendix A. Thus by Lemma \ref{lemma: cond argmax} we can conclude:
\begin{align*}
&\Pi_n\left( \sqrt{\frac{n}{\log n}}\left(\hat{V}_{\scriptscriptstyle{G}}(x)-V_0(x)\right) -Z^{\scriptscriptstyle{G_0}}_n(1) \leq a \mid Z_1, \ldots, Z_n \right)  \\
& \quad \stackrel{G_0}{\longrightarrow} \pb \left\{ \underset{t \in I_x}{\arg \max }\left\{t X-|t|^{\gamma+ 1}  K  -a t \right\} \leq 0 \right\} = \pb\{ X \leq a\},
\end{align*}
because $t \mapsto t(X-a) - |t|^{\gamma +1}K$ attains its maximum on $(-\infty,0]$ iff $X-a\leq 0$.
\end{proof}

We conclude this section by showing the implications of the above results for the uncertainty quantification of the posterior. Under the conditions of Theorem \ref{thm: BvM} (c.f.\ \cite{27}) 
\begin{align*}
	\frac{\sqrt{n}}{\sqrt{\log n}}\left(\hat{V}_n(x) - V_0(x)\right) \rightsquigarrow N\left(0, \frac{g_0(x)}{2 \gamma}\right).
\end{align*}
Thus as consequence of \eqref{eq: solution range equation} and the analogous definition of $\hat{F}_n(x)$, we have that (see proof of Corollary \ref{thm: BvM for F} in the Online Supplement \cite{64}):
\begin{align}\label{eq: convergence IIE for F}
	\frac{\sqrt{n}}{\sqrt{\log n}}\left(\hat{F}_n(x) - F_0(x)\right) \rightsquigarrow N\left(0, \frac{2 x}{\pi^2 \gamma} g_0(x)\right).
\end{align}
Together with \eqref{eq: convergence for F} below, this will imply that the credible intervals constructed from the Isotonized Inverse Posterior (IIP) of $\hat{F}_{\scriptscriptstyle{G}}(x)$ are (asymptotic) frequentist confidence sets. The following Corollary \ref{thm: BvM for F}, whose proof is given in the Online Supplement \cite{64}, summarizes formally the above discussion. Note that because the uncertainty quantification based on the IIP is always preferable to the one based on the Naive Bayes Posterior, we only propose uncertainty quantification using the latter.

\begin{corollary}[BvM for $\hat{F}_{\scriptscriptstyle{G}}$ and $F_{\scriptscriptstyle{G}}$ and UQ for $F_0$]\label{thm: BvM for F} Let $F_{\scriptscriptstyle{G}}$ as in \eqref{eq: naive post F0} and $\hat{F}_{\scriptscriptstyle{G}}$ as in \eqref{eq: IIP F0} and analogously let $F_n(x)$, with $V_n$ instead of $V_{\scriptscriptstyle{G}}$ in \eqref{eq: naive post F0} and $\hat{F}_n(x)$, with $\hat{V}_n$ instead of $\hat{V}_{\scriptscriptstyle{G}}$ in \eqref{eq: IIP F0}. Under the conditions of Theorem \ref{thm: BvM naive}, in probability under $G_0$ as $n \rightarrow \infty$:
\begin{align}
	&\frac{\sqrt{n}}{\sqrt{\log n}}\left(F_{\scriptscriptstyle{G}}(x)-F_n(x)\right) \mid Z_1, \ldots, Z_n \rightsquigarrow N\left(0, \frac{4 x}{\pi^2} g_0(x)\right). \numberthis \label{eq: convergence for naive F} 
\end{align}
Under the same conditions of Theorem \ref{thm: BvM}, in probability under $G_0$ as $n \rightarrow \infty$:
\begin{align}
    &\frac{\sqrt{n}}{\sqrt{\log n}}\left(\hat{F}_{\scriptscriptstyle{G}}(x)-\hat{F}_n(x)\right) \mid Z_1, \ldots, Z_n \rightsquigarrow N\left(0, \frac{2 x}{\pi^2 \gamma} g_0(x)\right).\numberthis \label{eq: convergence for F}
\end{align}
Further, if we define the posterior quantiles for isotonized posterior $\hat{F}_{\scriptscriptstyle{G}}$, for $\alpha \in [0,1]$, as:
\begin{align*}
	q_{\alpha} := \inf \{ t \in \mathbb{R} \, : \, \Pi_n (\hat{F}_{\scriptscriptstyle{G}}(x) \leq t \mid Z_1, \ldots, Z_n) \geq \alpha\}
\end{align*}	
then under $G_0$ as $n \rightarrow \infty$:
\begin{align}\label{eq: UQ}
    \pb \left( F_0(x) \in [q_{\alpha/2}, q_{1- \alpha/2}]\right) \rightarrow 1- \alpha.
\end{align}
\end{corollary}

\section{Simulation Study and Discussion about Uncertainty Quantification}

In this section, we present simulations to illustrate the practical behavior of the IIP in comparison to the IIE studied in \cite{27}. Furthermore we present a simulation study that compares the IIP-uncertainty quantification to the one proposed in \cite{19}.
We consider three distinct settings. The first two are \textit{synthetic} settings where the true underlying distribution function \( F_0 \) is known. The first two simulation settings serve as valuable benchmarks for assessing the quality of our methodology. In particular, the first one verifies that the variance attained by the IIP coincides with the theoretical one that we derived in Theorem \ref{thm: BvM for F}. The second \textit{synthetic} setting compares the uncertainty quantification of our methodology to the bootstrap-based one derived in \cite{19}. The last simulation setting is a \textit{real-data} application involving the positions of stars in the Globular Cluster Messier 62. In the last case, the ground truth distribution of the stars is unknown. Also in this setting, we compare our methodology with the bootstrap of \cite{19}. 

\begin{figure}[H]
	\centering
    \begin{subfigure}[t]{0.48\linewidth}
        \centering
        \includegraphics[width=\linewidth]{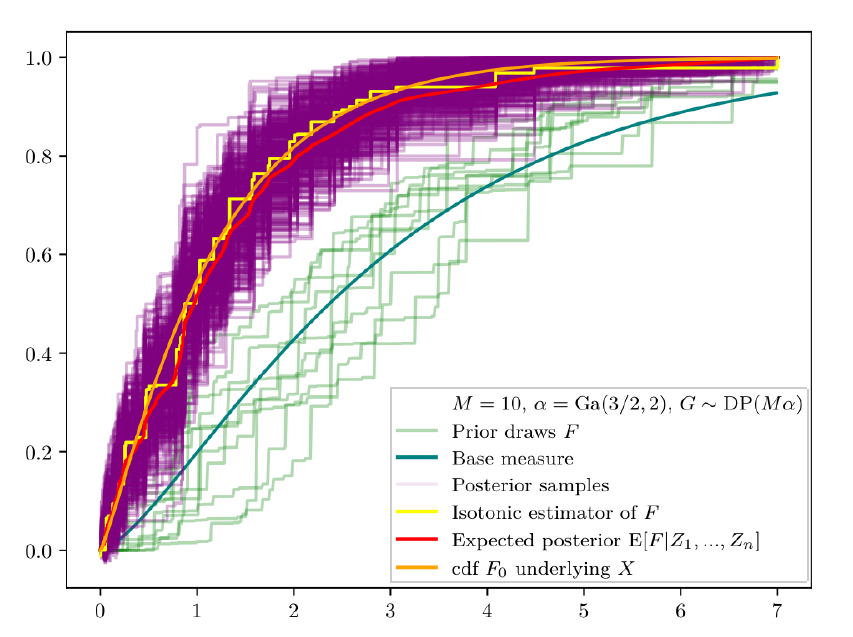}
        \caption*{$n=200$}
    \end{subfigure}
    \hspace{-0.2cm}
    \begin{subfigure}[t]{0.48\linewidth}
        \centering
        \includegraphics[width=\linewidth]{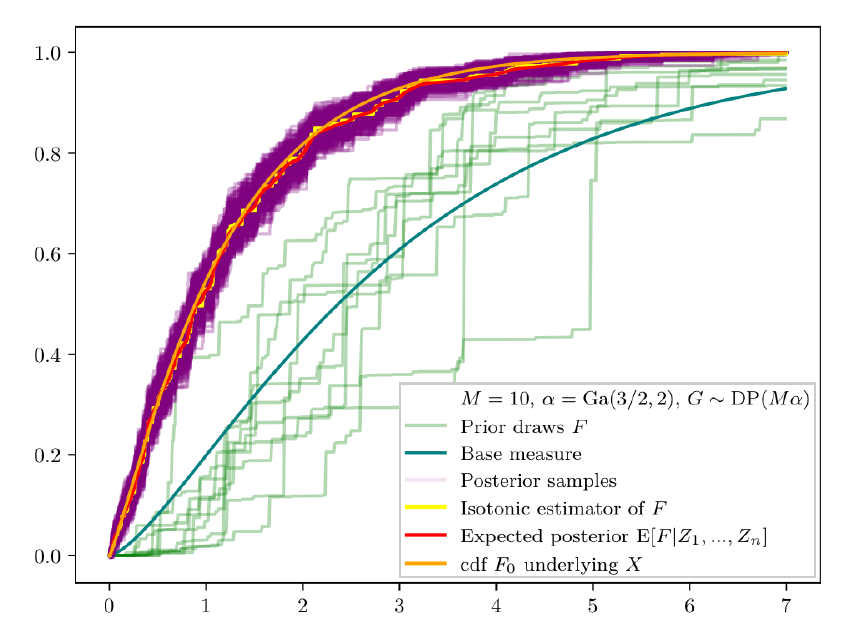}
        \caption*{$n=2000$}
    \end{subfigure}
	\\
    \centering
	\begin{subfigure}[t]{0.5\linewidth}
        \centering
        \includegraphics[width=\linewidth]{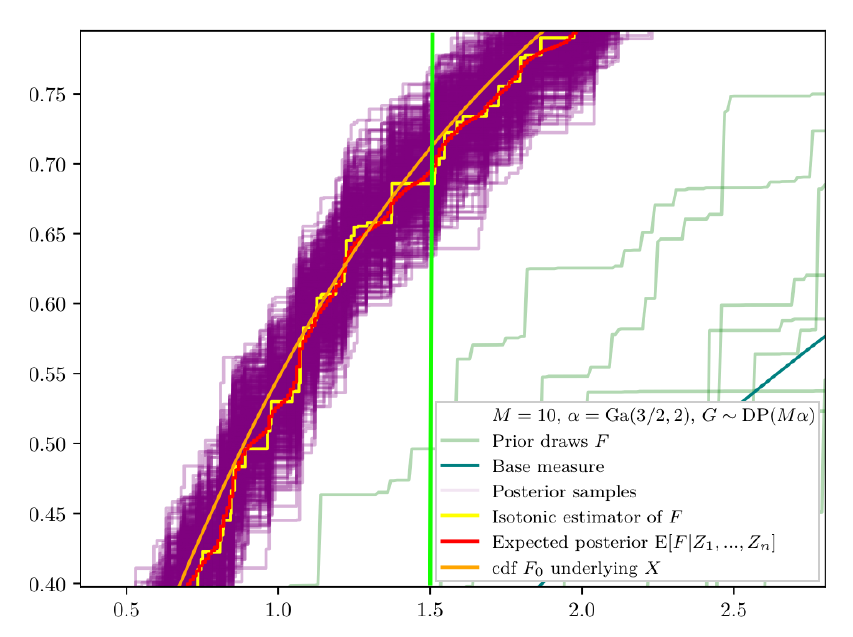}
		\caption*{Zoom for $n=2000$}
    \end{subfigure}
	\\
	\vspace{0.2cm}
    \begin{subfigure}[t]{0.48\linewidth}
        \centering
        \includegraphics[width=\linewidth]{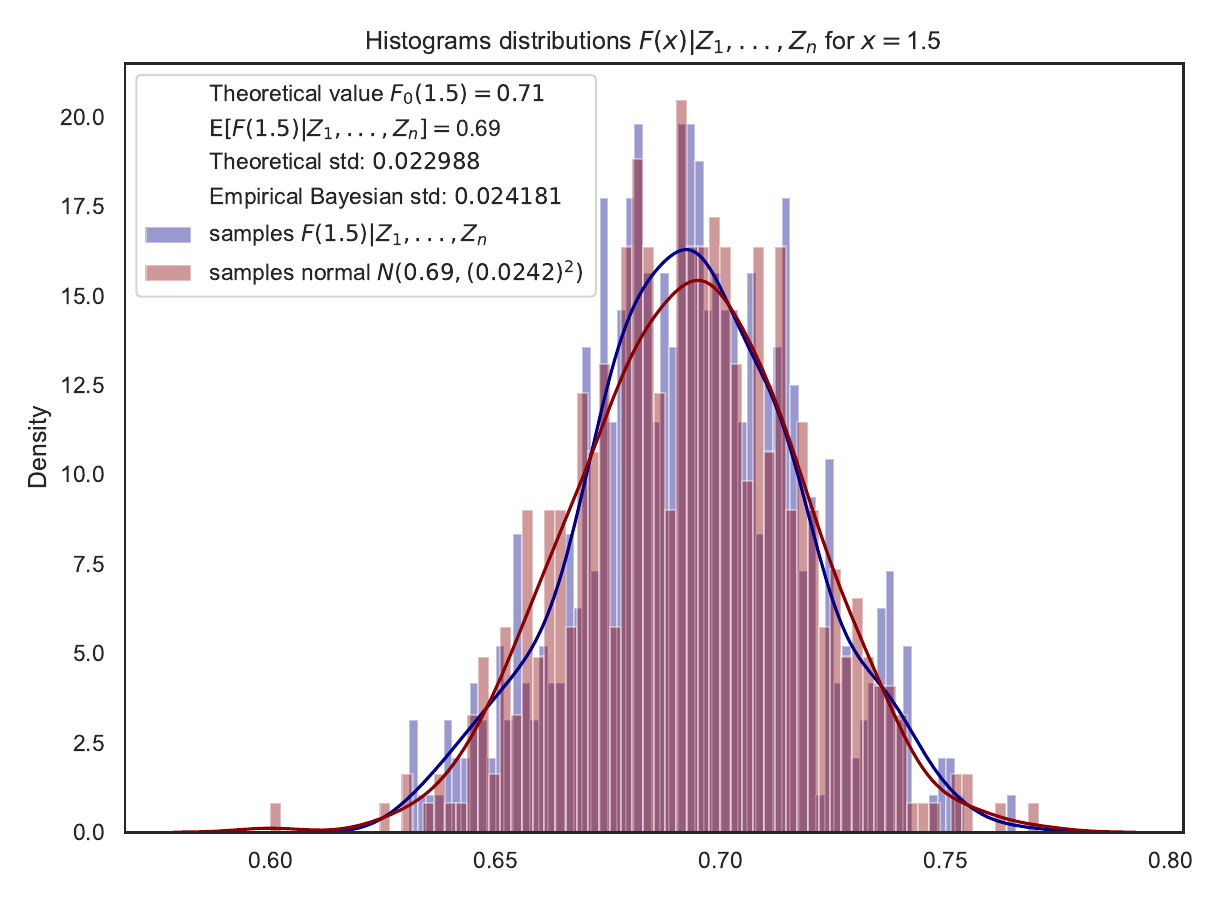}
        \caption*{Histograms}
    \end{subfigure}
    \hspace{-0.3cm}
    \begin{subfigure}[t]{0.52\linewidth}
        \centering
		\vspace*{-5.3cm}
        \includegraphics[width=\linewidth]{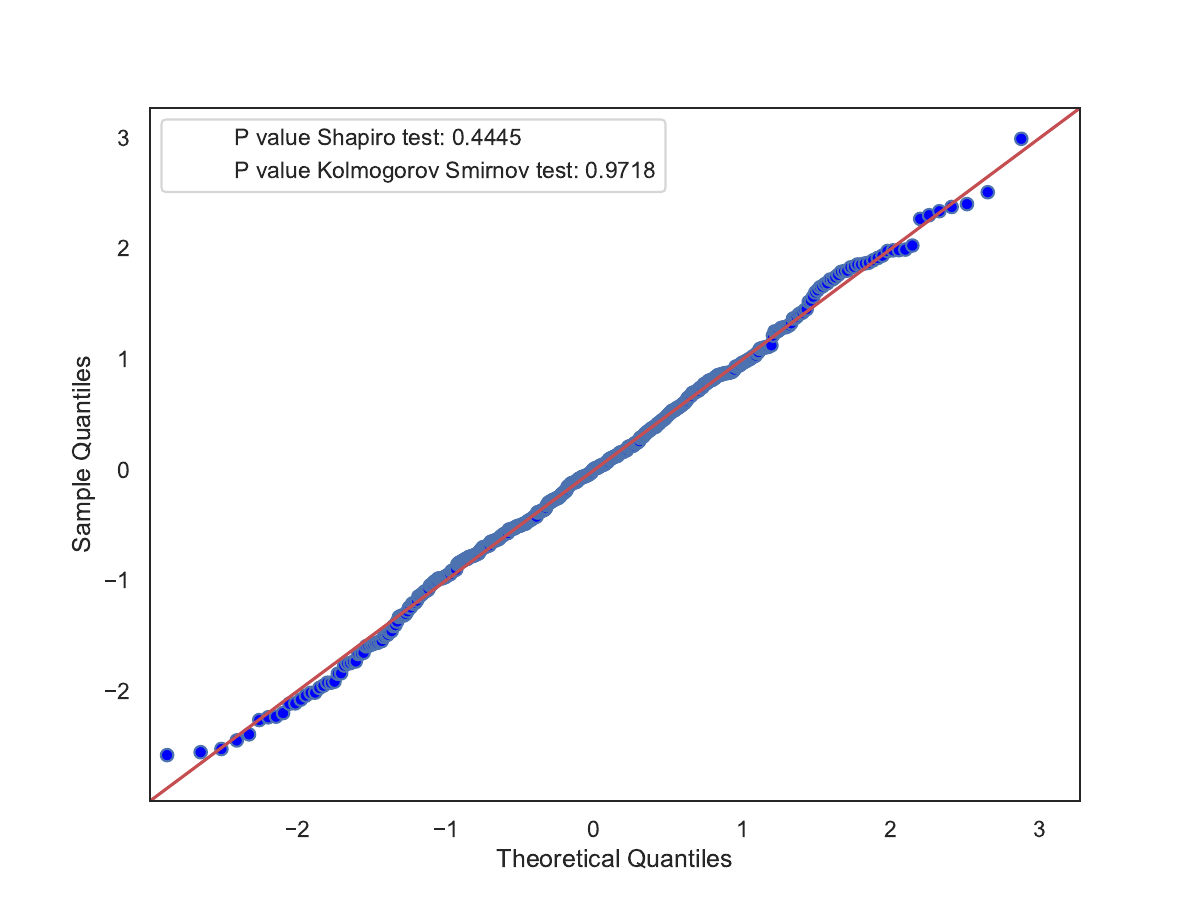}
        \caption*{QQ-plot}
    \end{subfigure}
	\caption{\textit{Synthetic} setting 1: uncertainty for different sample sizes and analysis of normality at $x=1.5$, based on 300 draws from the posterior.}\label{fig: analysis normality}
\end{figure}

In Figure \ref{fig: analysis normality}, we analyze the setting where the underlying cdf \( F_0 \) is exponential with rate 1.2. Samples \( Z_1, \ldots, Z_n \) used to compute the projection posterior of the IIP are generated hierarchically: \( X_1, \ldots, X_n \sim F_0 \), \( Y_1, \ldots, Y_n \sim \mathrm{Be}(1, 1/2) \), with \( Z_i = X_i \cdot Y_i \) for all \( i \). The plots use the following color scheme: green for estimates of \( F \) based on prior draws from \( G \sim \operatorname{DP}(M \alpha) \) (parameters \( M \) and \( \alpha \) specified in the plots); purple for draws from the isotonized posterior; red for the average of these draws (approximating the posterior mean);  yellow for the Isotonic Inverse Estimator (IIE) of \cite{27}; and orange for \( F_0 \).  
\vspace{-0.5cm} 
The first row of Figure \ref{fig: analysis normality} compares the behavior of the IIP and the IIE as the sample size goes from 200 to 2000. With \( n = 200 \), the IIE and the average of the IIP draws already behave very similarly and with \( n = 2000 \), the uncertainty significantly reduces and the two estimation procedures essentially coincide. The average IIP appears to be somewhat smoother, which will typically be seen as an advantage. The middle row focuses on \( n = 2000 \), providing a zoomed-in view and highlighting \( x = 1.5 \), the point used for subsequent analysis of the normality. The final row examines normality (\( n = 2000 \)): histograms, QQ-plots, and statistical tests strongly confirm normality of the gathered samples. Furthermore, the Bayesian-estimated standard deviation aligns closely with the theoretical value attained by the frequentist procedure (c.f.\ \cite{27}, note that \( \gamma = 1 \) because \( F_0 \) is Lipschitz). In all the plots we employ 300 posterior draws. 
\vspace{-0.5cm} 
\begin{figure}[H]
	\centering
    \begin{subfigure}[t]{0.52\linewidth}
        \centering
        \includegraphics[width=\linewidth]{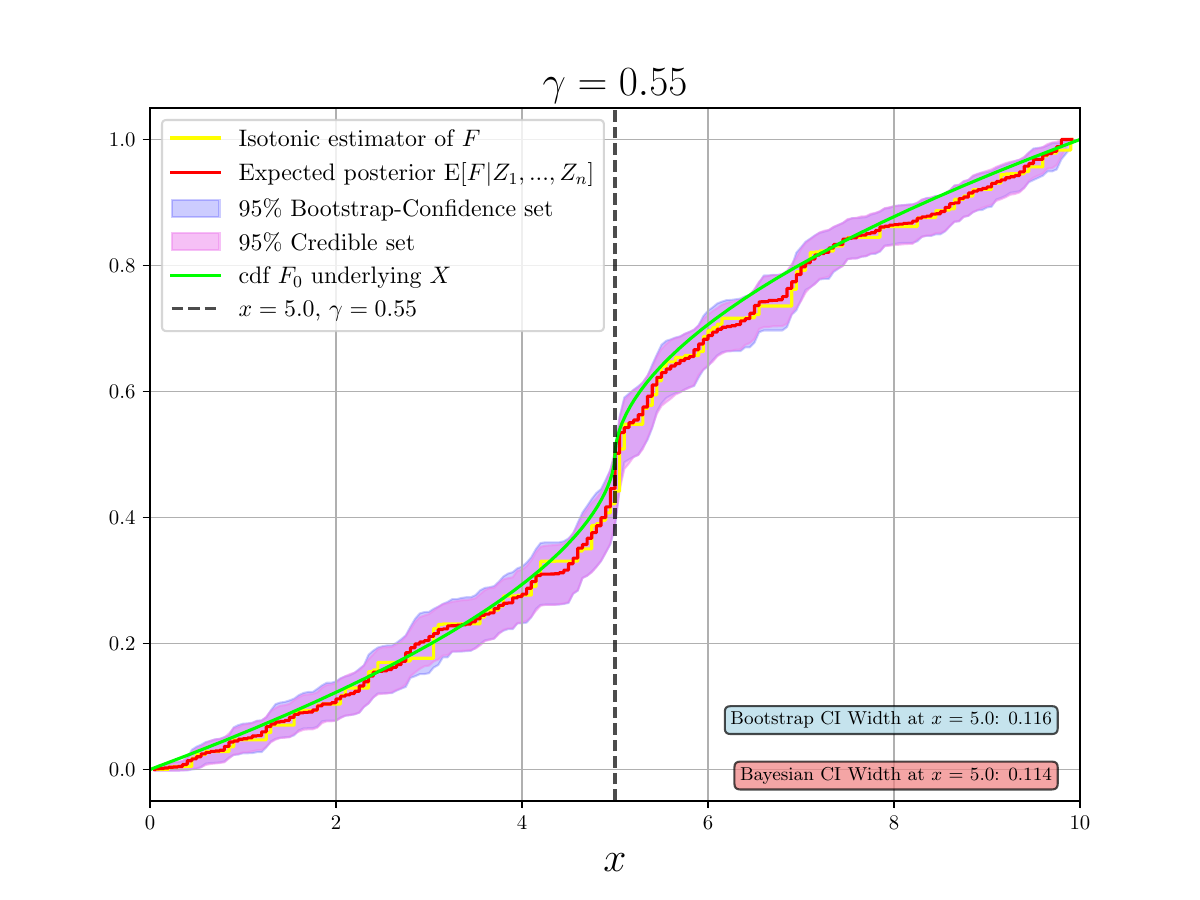}
    \end{subfigure}
    \hspace{-0.8cm}
    \begin{subfigure}[t]{0.52\linewidth}
        \centering
        \includegraphics[width=\linewidth]{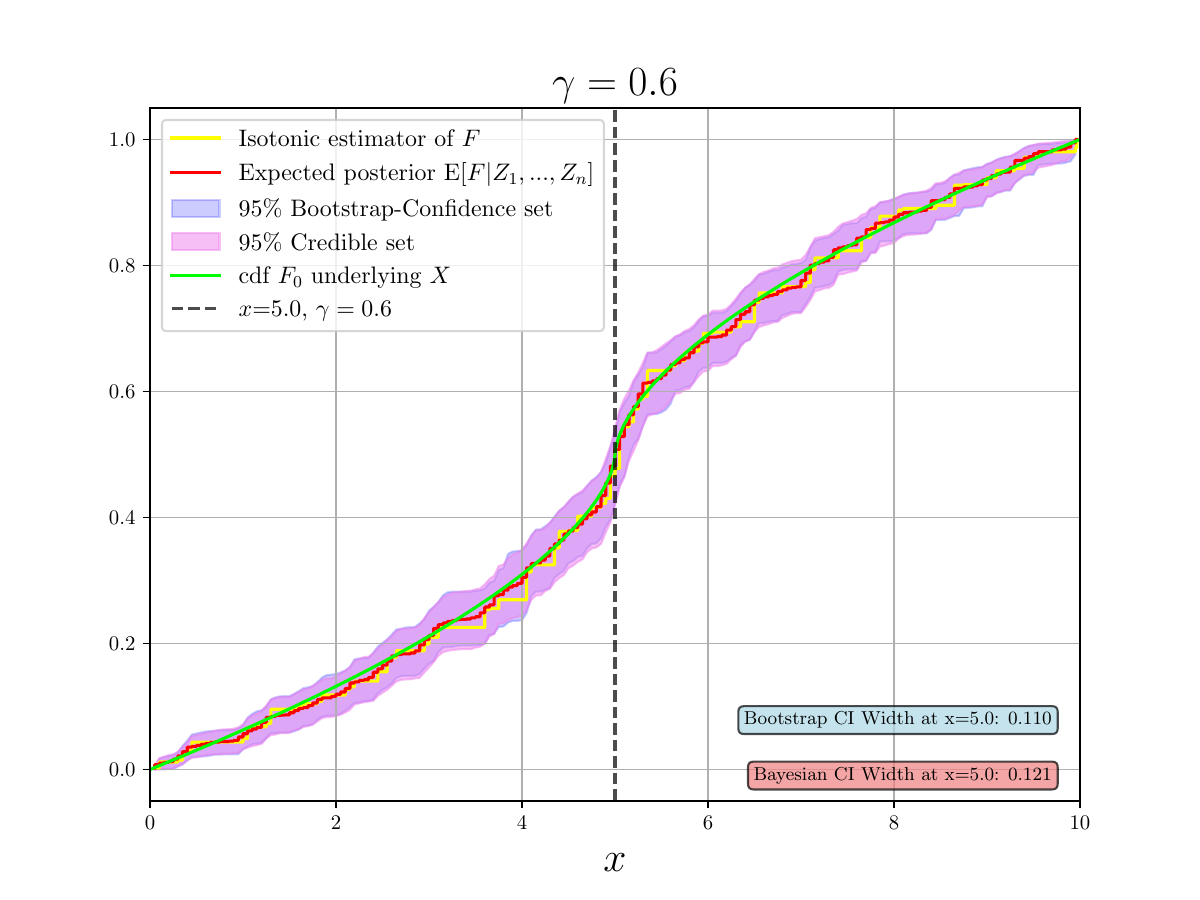}
    \end{subfigure}
	\\
    \begin{subfigure}[t]{0.52\linewidth}
        \centering
        \includegraphics[width=\linewidth]{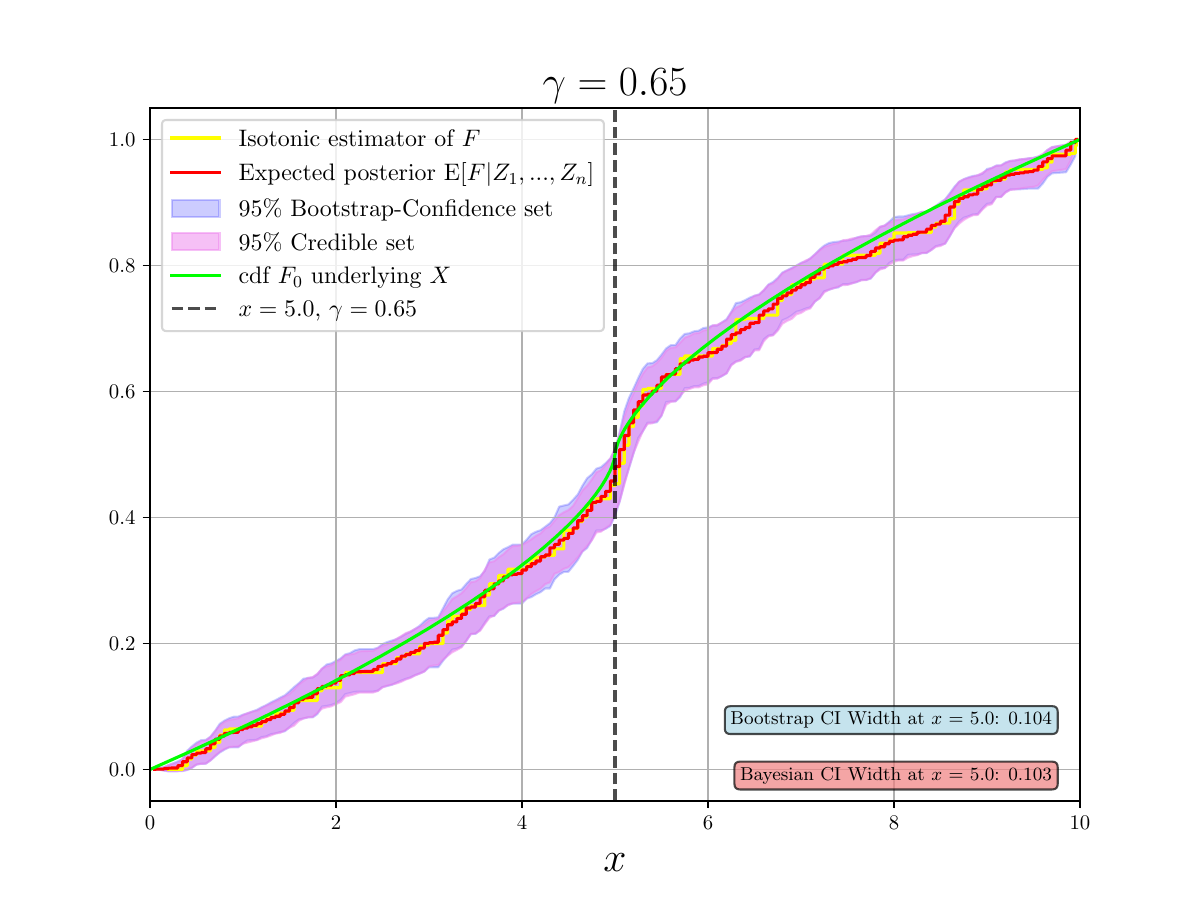}
    \end{subfigure}
    \hspace{-0.8cm}
    \begin{subfigure}[t]{0.52\linewidth}
        \centering
        \includegraphics[width=\linewidth]{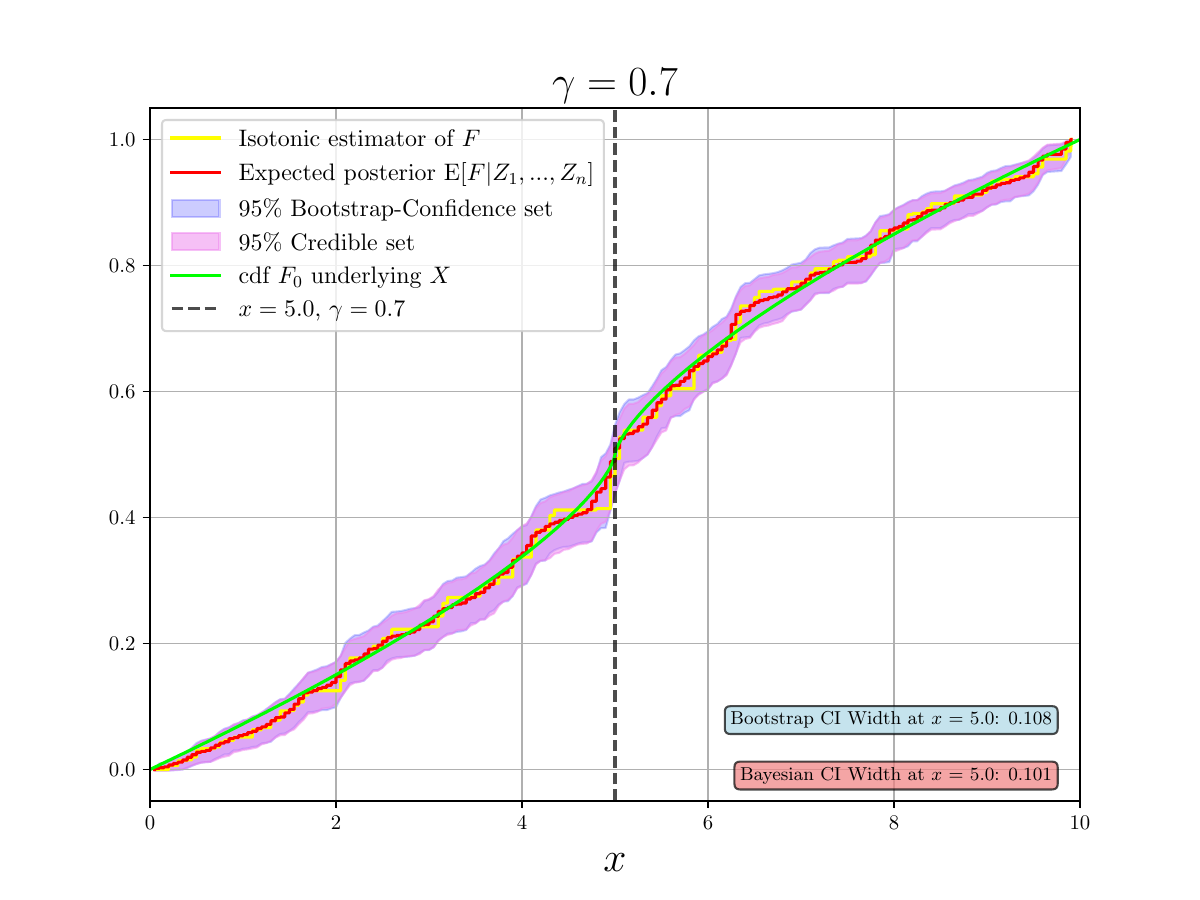}
    \end{subfigure}
	\\
    \begin{subfigure}[t]{0.52\linewidth}
        \centering
        \includegraphics[width=\linewidth]{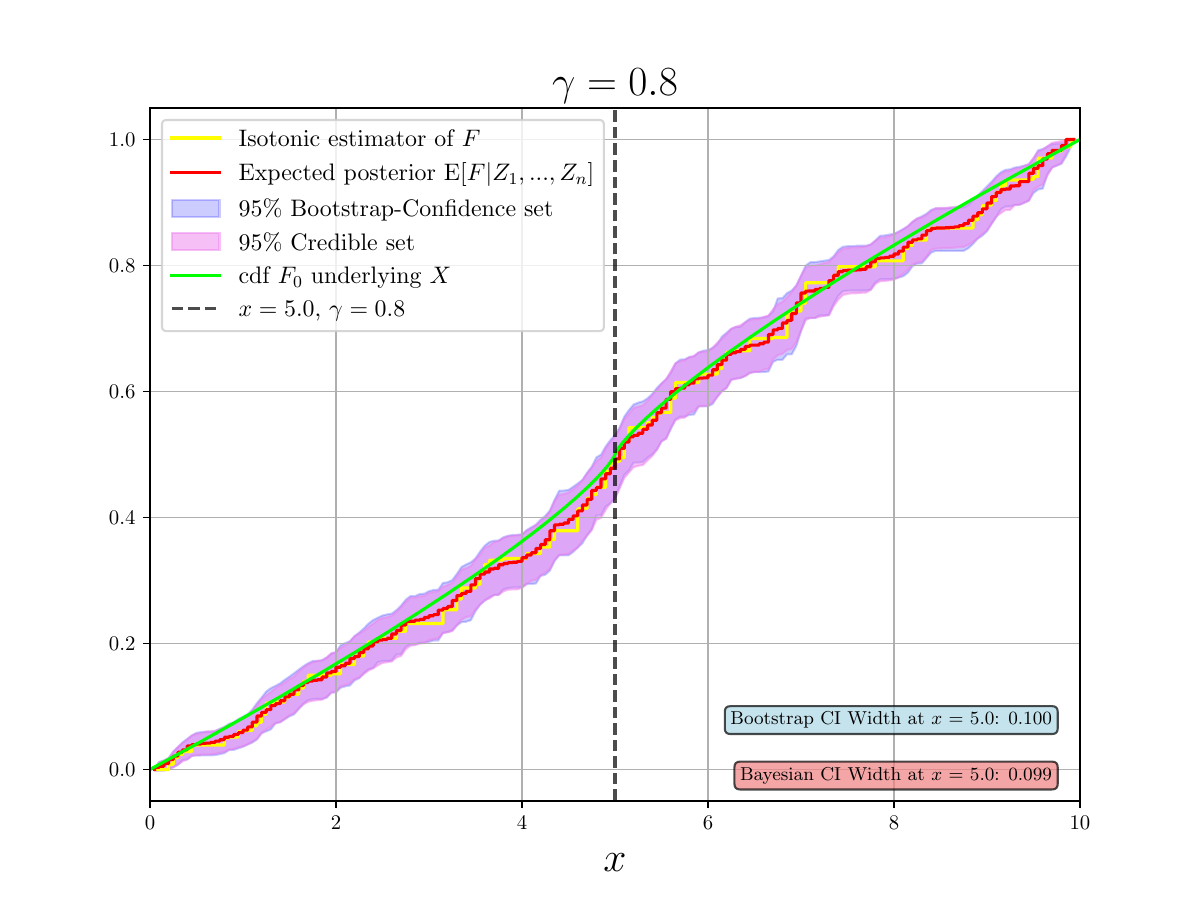}
    \end{subfigure}
    \hspace{-0.8cm}
    \begin{subfigure}[t]{0.52\linewidth}
        \centering
        \includegraphics[width=\linewidth]{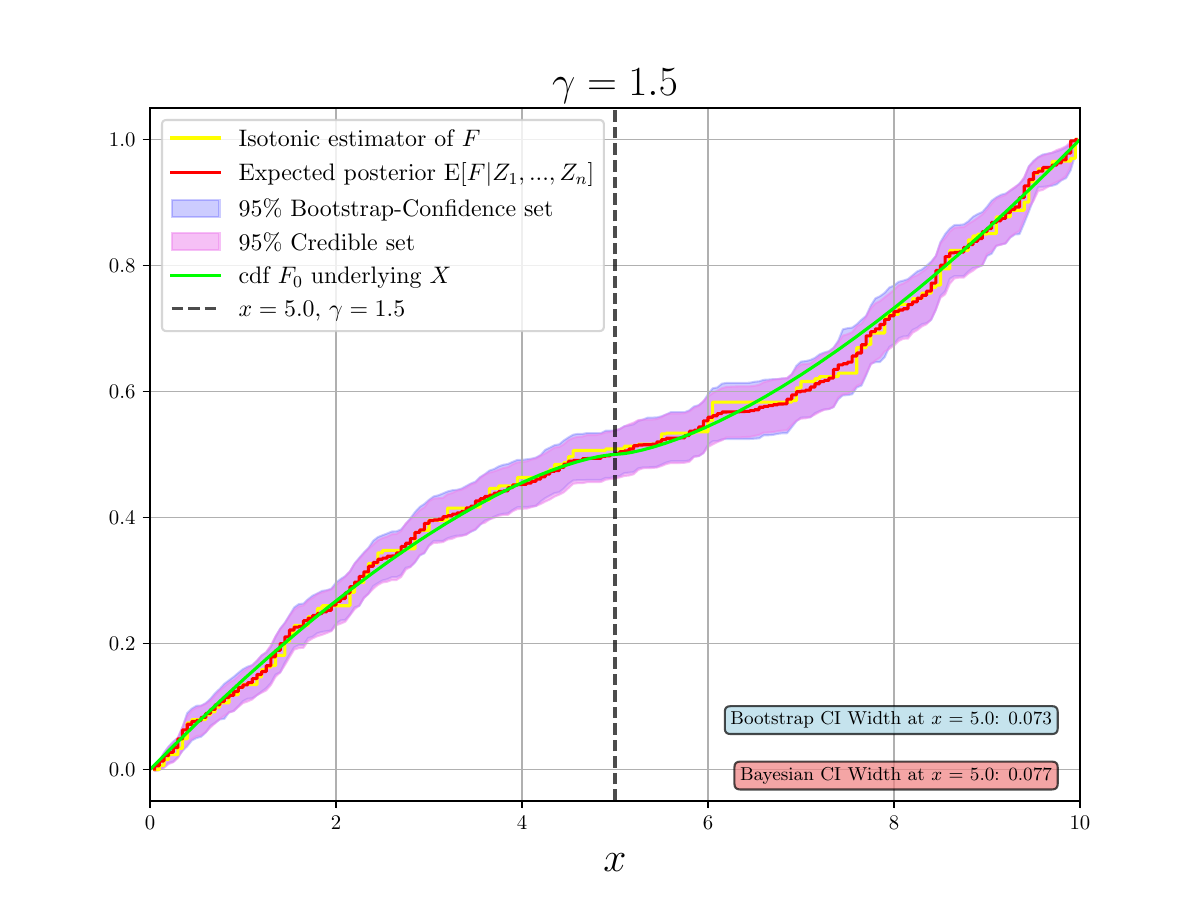}
    \end{subfigure}
	\caption{\textit{Synthetic} setting 2: comparison uncertainty quantification IIP and bootstrap method from \cite{19}.}\label{fig: compare UQ}
\end{figure}
\vspace{-1cm}
\begin{table}[H]
\centering
\begin{tabular}{c|cccccc}
   & $\gamma = 0.55$ & $\gamma = 0.6$ & $\gamma = 0.65$ & $\gamma = 0.7$ & $\gamma = 0.8$ & $\gamma = 1.5$ \\ \hline
Bootstrap CI width at $x = 5$ &  0.116   &  0.110    &  0.104   &  0.108   &   0.100  &   0.073  \\
Bayesian CI width at $x = 5$ &   0.114 &  0.121   &   0.103  &   0.101  &   0.099  &   0.077  \\
\end{tabular}
\label{tab:your_label}
\end{table}

In Figure \ref{fig: compare UQ}, we analyze six different settings where at $x=5$ the underlying cdfs $F_0$ have varying degrees of smoothness, making them locally Hölder smooth of degree $\gamma \in \{ 0.55, 0.6, 0.65, 0.7, 0.8, 1.5 \}$. In particular, for an appropriate constant $K>0$, $F_0(y) = \frac{1}{2} - K (5-y)^\gamma$, for $y \in [0,5]$ and $F_0(y) = \frac{1}{2} + K (y-5)^\gamma$, for $y \in [5,10]$. The plots use the following color scheme: pink for the Bayesian-based credible sets; red for the approximated posterior mean; light blue for the bootstrap-based confidence sets of \cite{19}; yellow for the Isotonic Inverse Estimator (IIE) of \cite{27}; and lime green for \( F_0 \). The prior base measure used for these simulation is the same used in the previous simulation study. In all cases \( n = 2000 \) and we employed 300 posterior-based and bootstrap-based draws. This comparison confirms that the IIP delivers valid uncertainty quantification for varying degrees of smoothness and that the uncertainty quantification delivered by \cite{19} behaves closely, even though the authors formally showed it works only for a differentiable underlying truth $F_0$.
\vspace{-0.5cm} 
\begin{figure}[H]
    \centering
	\begin{subfigure}{0.55\linewidth}
        \includegraphics[width=\linewidth]{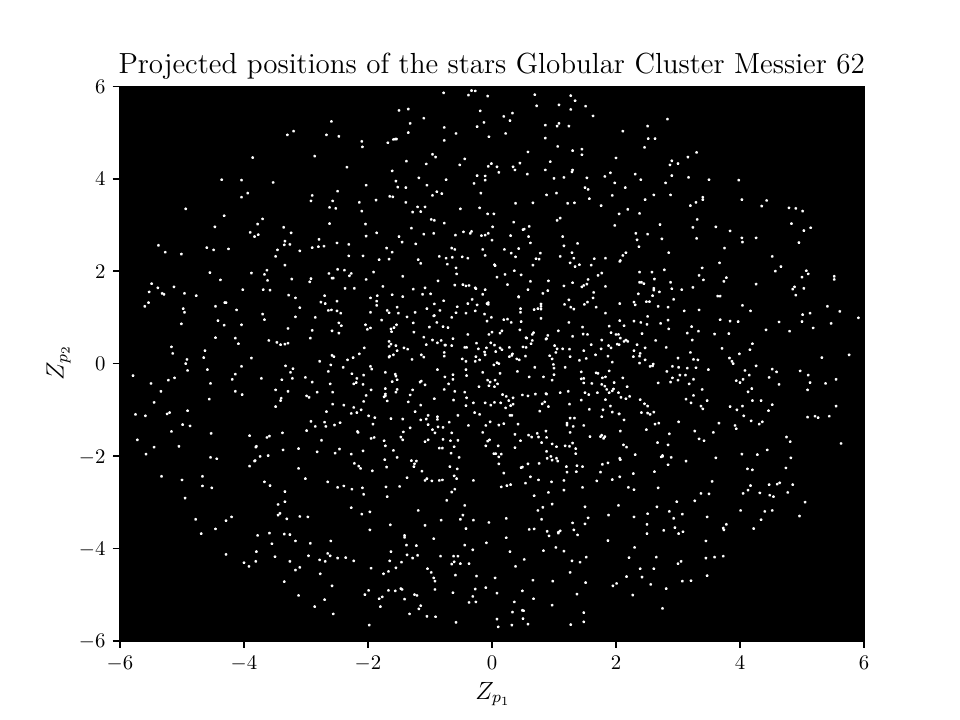}
    \end{subfigure}

    \begin{subfigure}{0.51\linewidth}
        \includegraphics[width=\linewidth]{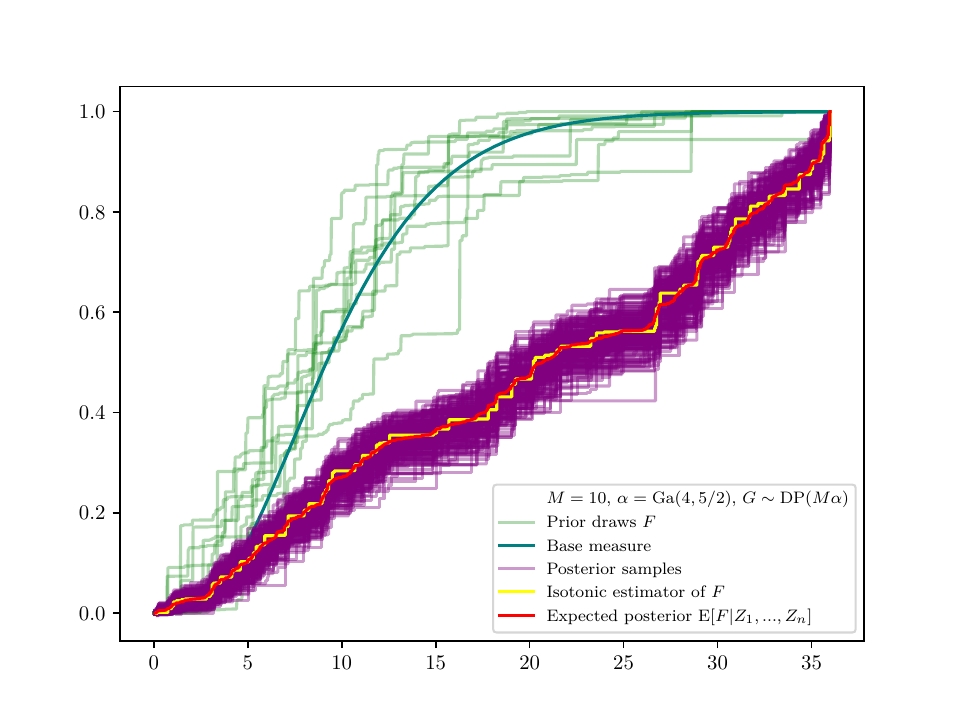}
    \end{subfigure}  
	\hspace{-0.6cm}
    \begin{subfigure}{0.515\linewidth}
        \includegraphics[width=\linewidth]{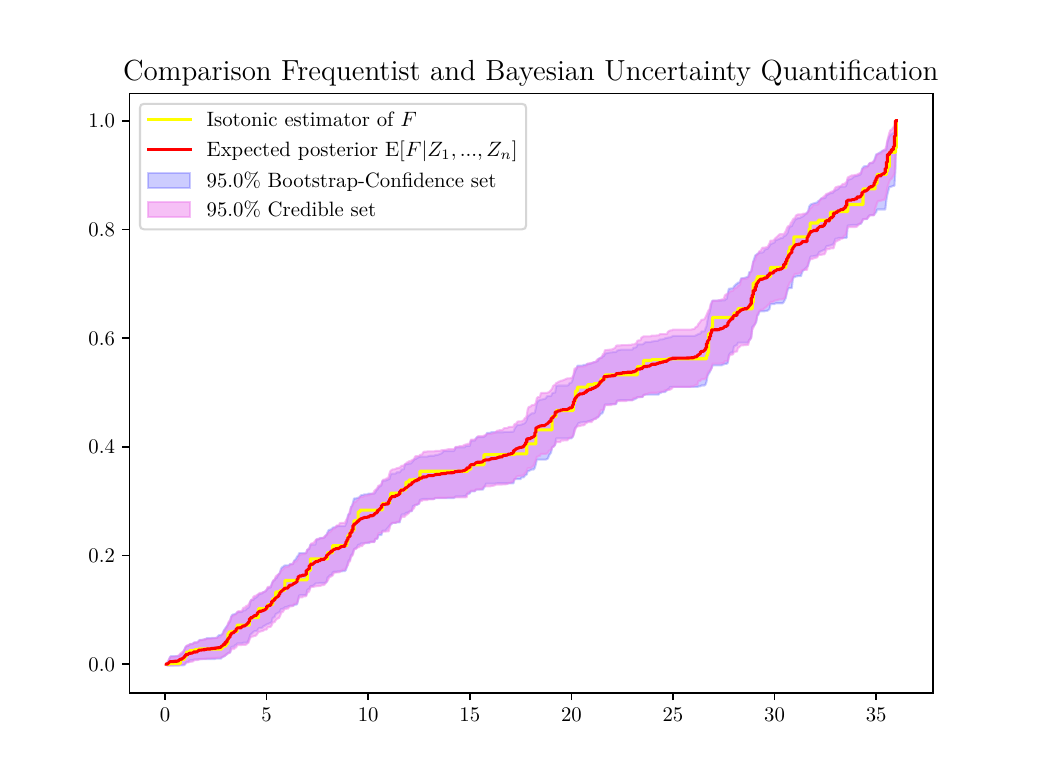}
    \end{subfigure}

    \caption{\textit{Real-data} application: $n=1400$ projected star positions. On the left-hand side, the IIE from \cite{27} and the IIP with 300 posterior samples; on the right-hand side, a comparison between the $95\%$ bootstrap confidence set for the IIE from \cite{19} and the $95\%$ credible sets for the IIP.}
    \label{fig:stars}
\end{figure}
\vspace{-0.5cm}
In Figure \ref{fig:stars}, we apply the methodology to star position data from globular cluster Messier 62. Let \( \mathbf{Z} = (Z_{p_1}, Z_{p_2}, Z_{p_3}) \) represent a star's three-dimensional position, where only \( (Z_{p_1}, Z_{p_2}) \) are observed. Here, \( X = Z_{p_1}^2 + Z_{p_2}^2 + Z_{p_3}^2 \sim F_0 \), and \( Z = Z_{p_1}^2 + Z_{p_2}^2 \). Beyond serving as a compelling application of our methodology, this example highlights a few key conclusions. First of all, also for real data the IIE closely resembles the IIP. However, because in this case the true $\gamma$ is unknown, it is not immediately clear how the frequentist uncertainty quantification compares with the Bayesian one. From the work of \cite{19}, we know that bootstrap procedures work for Wicksell's problem, at least when the underlying truth is differentiable at $x$. We thus compared our methodology with the one proposed in \cite{19} and we observed that they behave very similarly with $1000$ or more samples. In Figure \ref{fig:stars} the reader can observe an instance with $n=1400$.

\subsubsection*{Discussion about Uncertainty Quantification}

The main limitation of the approach in \cite{27} is the lack of uncertainty quantification, which is essential for practical applications. The asymptotic variance of the limiting distribution achieved by the IIE involves two nuisance parameters: the local smoothness $\gamma$ and the density $g_0(x)$ of the observed data at the evaluation point $x$. While $g_0(x)$ may be estimable (for most practical applications the rate of convergence would be slower than $\sqrt{n}/\sqrt{\log n}$), estimating $\gamma$ is particularly challenging, as it requires recovering the local smoothness of the target distribution function $F_0$ at a specific point $x$ using samples from $g_0$—effectively a separate inverse problem. To our knowledge, no practical methods currently exist for this task. Moreover, because of the method proposed in this paper, estimating $\gamma$ offers little practical benefit and thus makes this research direction of limited value. This highlights the key advantage of the Bayesian approach proposed in this paper: it bypasses the need to estimate either $\gamma$ or $g_0(x)$, as the posterior distribution inherently provides valid uncertainty quantification over a wide range of local smoothness levels, i.e.\ all $\gamma > 1/2$. For $\gamma \leq 1/2$, the inverse problem becomes severely ill-behaved due to the unboundedness of the density of the observations $g_0$. This might even lead to a slower rate of convergence both for the naive plug-in estimator and for the IIE. Therefore, for all practically relevant cases with $\gamma > 1/2$, practitioners can apply our method without needing to estimate nuisance parameters or make strong assumptions about the form of $F_0$, thanks to its adaptive nature and the automatic uncertainty quantification.

The results displayed in Figure \ref{fig: compare UQ} suggest that the bootstrap method of \cite{19} behaves very closely to our method. We can thus make the empirically grounded hypothesis that the bootstrap method proposed in \cite{19} might work under the less stringent assumptions of the current paper. The value of this simulation study has to be taken carefully, as a formal mathematical proof is still needed. However, proving this hypothesis is nontrivial and we therefore leave this as potential future work. 

From a practitioner's perspective, thus, the key advantage of our Bayesian methodology lies in its provable robustness across a broad class of smoothness levels: it provides valid uncertainty quantification for all $\gamma > 1/2$. In contrast, the bootstrap approach in \cite{19} is provably consistent only when $F_0$ is differentiable at $x$, a substantially stronger condition (stricter than $\gamma = 1$). This distinction is crucial, since determining whether $F_0$ is differentiable at a given point based solely on samples from $g_0$ constitutes a challenging inverse problem—both theoretically and in practice. By avoiding this requirement, our method offers a more broadly applicable and practitioner-friendly solution.

\section{Conclusion}

We have introduced a novel Bayesian approach for nonparametric estimation in Wicksell’s problem, leveraging a Dirichlet Process prior on the distribution of observables rather than unobservables. This deviation from the classical Bayesian framework in inverse problems allows for computational efficiency through posterior conjugacy. By projecting the posterior draws onto the \( \mathbb{L}_2 \) subspace of increasing, right-continuous functions, we obtained the Isotonized Inverse Bayes Posterior (IIP) and we verified that it allows for a Bernstein–von Mises theorem with minimax asymptotic variance, reflecting the Hölder continuity of the true distribution. This proves that our methodology gives asymptotically correct uncertainty quantification under the frequentist point of view and it highlights its potential for Bayesian nonparametrics in inverse problems. 

Regarding potential future work, several interesting directions emerge. First, it would be important to verify whether the bootstrap procedure of \cite{19} remains valid for all degrees of smoothness of the underlying truth, not just for $F_0$ differentiable at $x$. Second, an open question is whether a Bernstein–von Mises (BvM) result can be obtained for Wicksell’s problem within the classical Bayesian framework—that is, by placing a Dirichlet Process (DP) prior $F$ on \( F_0 \). Formally obtaining a nonparametric hierarchical Bayesian model for the observations $Z_1,\ldots,Z_n$ as follows:
\begin{itemize}
    \item[(i)] A probability distribution $F$ on $\mathbb{R}_+$ is generated from the Dirichlet Process prior $\operatorname{DP}(\alpha)$.
    \item[(ii)] An i.i.d.\ sample $X_1,\ldots,X_n$ is generated from $F$.
    \item[(iii)] An i.i.d.\ sample $Y_1,\ldots,Y_n$ is generated from a density of a beta distribution $\operatorname{Be}(1,1/2)$, independent of the other samples $X_1,\ldots,X_n$.
    \item[(iv)] The observations are $Z_i = Y_i \cdot X_i$, for $i=1,\ldots,n$, with density $g_{\scriptscriptstyle{F}}(z) := \int_z^{\infty} \frac{dF(x)}{2\sqrt{x^2 - xz}}$.
\end{itemize}
In the above model, the posterior is a Mixture of Dirichlet Processes (MDP), with distribution:
\begin{align*}
    &F \mid Z_1,\ldots,Z_n \sim \int \mathrm{DP}\bigg( \alpha + \sum_{i=1}^{n} \delta_{X_i}\bigg) \: d \mathrm{P}(X_1,\ldots,X_n \mid Z_1,\ldots,Z_n), \numberthis \label{eq: posterior distribution}
\end{align*}
where \eqref{eq: posterior distribution} follows directly from Corollary 3.1 in \cite{21}. The authors of the present paper have obtained unpublished preliminary results in this model, showing:
\begin{thmlist}
	\item denoting by $h$ Hellinger distance, there exists a constant $M$ such that, in $G_0$-probability: \label{thm: contraction rate hellinger}
	$$
		\Pi_n\Big(F: h\big(g_{\scriptscriptstyle{F}}, g_0 \big) \geq M n^{-1 / 4} (\log n)^{1/4} \mid Z_1, \ldots, Z_n\Big) \stackrel{G_0}{\longrightarrow} 0,
	$$
	\item denoting by $V_{\scriptscriptstyle{F}}$ the posterior of $V$ with the above specified prior, for every $x \geq 0$, there exists a constant $M$ such that, in $G_0$-probability:\label{thm: contraction rate V}
	\begin{align*}
		\Pi_n\left(F : |V_{\scriptscriptstyle{F}}(x) - V_0(x)| \geq M n^{-1 / 4} (\log n)^{3/4} \mid Z_1, \ldots, Z_n\right) \stackrel{G_0}{\longrightarrow} 0.
	\end{align*}
\end{thmlist}
The contraction rate $n^{-1/4}$ for the Hellinger distance is equal to the $\mathbb{L}_2$-minimax rate for estimating a density in a Sobolev space of order $1/2$ (see \cite{45}). Because the densities $g_0$ are in $H^{1/2}[0,\infty)$ if $F_0$ possesses a bounded density, this suggests that the rate  $n^{-1/4}$ is close to optimal. The pointwise contraction rate for the functional \( V_0 \) follows from the contraction rate in Hellinger distance, combined with a proof technique based on a Laplace tranform rewriting of the posterior as in eq.\ (2.13) in \cite{37}, a local perturbation and a LAN expansion as in \cite{27}. This perturbation, which is believed to correspond to the efficient influence function for the functional of interest, is essential for handling arguments involving contraction rates. However, the local nature of the perturbation in Wicksell's problem comes at the cost of contracting at the Hellinger distance rate $h$, which is too slow (and not $h^2 = O_p (\sqrt{\log{n}} / \sqrt{n})$ as in the proof style of \cite{37}). Even though the authors have made significant efforts to overcome this issue, this challenge remains unsolved. 

We conclude by remarking that the work of the present paper introduces a fresh perspective on Bayesian inverse problems by employing the concept of the \textit{projection posterior}, first proposed in \cite{54}. The methods proposed here are also significantly computationally faster if compared to the standard hierarchical Bayes model illustrated above, because they do not rely on computationally expensive Gibbs samplers. Finally, a valuable direction for future research is to identify other problems where this approach could be effectively utilized.

\begin{appendix}
\section{Complementary results}
\textit{Measure-theoretical framework}: in this paper, the sample space is $(\mathbb{R}^+, \mathscr{B}(\mathbb{R}^+))$, and $\mathfrak{M}$ (with Borel $\sigma$-field $\mathscr{M}$ for the weak topology) is the collection of all Borel probability measures on $(\mathbb{R}^+, \mathscr{B}(\mathbb{R}^+))$. The Dirichlet process prior is a probability measure on $\mathfrak{M}$. A prior on the set of probability measures $\mathfrak{M}$ is a probability measure on a sigma-field $\mathscr{M}$ of subsets of $\mathfrak{M}$. Alternatively, it can be viewed as map from some probability space $(\Omega, \mathscr{U}, \mathrm{Pr})$ into $(\mathfrak{M}, \mathscr{M})$. We can think of the hierarchy of the spaces as follows: $(\Omega, \mathscr{U}, \mathrm{Pr}) \rightarrow (\mathfrak{M}, \mathscr{M}) \rightarrow (\mathbb{R}^+, \mathscr{B}(\mathbb{R}^+))$. We put a Dirichlet process prior over $G$ and write $G \sim \Pi = \mathrm{DP}( \alpha)$ for some base measure $\alpha$ and prior precision $|\alpha|$. This means that for any measurable set $B \in \mathscr{M}$ we have $\Pi(B) = \mathrm{Pr} ( G \in B)$. Throughout the paper, we use the notations $\pb$ and $\ex$, which should be understood as the outer probability measure and outer expectation, respectively, whenever issues of measurability in the considered processes arise (see \cite{5} for a rigorous treatment of these topics). 

Furthermore, we use the notion of conditional convergence in distribution. This happens if the bounded Lipschitz distance between the conditional law of Borel measurable $X_n$ given $\mathcal{B}_n$ and a Borel probability measure $L$ tends to zero, where the convergence can be in outer probability. The sequence $X_n$ tends to $L$ conditionally given $\mathcal{B}_n$ in outer probability if
\begin{align}\label{eq: formal definition conditional convergence}
d_{\scriptscriptstyle{ B L}}\left(\mathcal{L}\left(X_n \mid \mathcal{B}_n\right), L\right) :=\sup _{f \in B L_1}\left| \, \mathrm{E}\left(f\left(X_n\right) \mid \mathcal{B}_n\right)-\int f d L\right| \xrightarrow{\mathrm{P}} 0,
\end{align}
where $ B L_1$ is the space of bounded Lipschitz functions with Lipschitz constant 1.

In particular if $Z_1,\ldots,Z_n \overset{\text{i.i.d.}}{\sim} \pb_{\scriptscriptstyle{Z}}$, then $X_n \mid Z_1,\ldots,Z_n \rightsquigarrow N(0,\sigma^2)$ means:
\begin{align}\label{eq: formal definition conditional convergence 1}
	d_{\scriptscriptstyle{ B L}}\left( \pb (X_n \in \cdot \mid Z_1,\ldots,Z_n), N(0,\sigma^2) \right) \xrightarrow{\pb_{\scriptscriptstyle{Z}}} 0.
\end{align}

\begin{lemma}\label{eq: lemma equiv to claim}
	\eqref{eq: equiv to claim} implies \eqref{eq: claim}. 
\end{lemma}
	\begin{proof}
	By Theorem 1.13.1 (i) and (ii) in \cite{5} the claim \eqref{eq: claim} is obtained as soon as we show (c.f.\ \eqref{eq: formal definition conditional convergence}), for $X \sim N\left(0, \frac{g_0(x)}{2 \gamma}\right)$:
	\begin{align*}
		\sup _{f \in B L_1}\left| \, \mathrm{E}\left(f\bigg(\frac{\sqrt{n}}{\sqrt{\log n}}\big(\hat{V}_{\scriptscriptstyle{G}}(x)-\hat{V}_n(x)\big)\bigg) \mid Z_1,\ldots,Z_n \right)- \ex \left(f(X) \right) \right| \xrightarrow{G_0} 0.	
	\end{align*}
	By Lemma 2.1 given in the Online Supplement \cite{64}, for any bounded 1-Lipschitz function $f : \mathbb{R} \mapsto \mathbb{R}$ and $Z^{\scriptscriptstyle{G_0}}_n$ as in \eqref{eq: Z_n process}:
	\begin{align*}
		&\bigg|\mathrm{E} \bigg[f\bigg(\frac{\sqrt{n}}{\sqrt{\log n}}\big(\hat{V}_{\scriptscriptstyle{G}}(x)-\hat{V}_n(x)\big)\bigg) - f\bigg(\frac{\sqrt{n}}{\sqrt{\log n}}\big(\hat{V}_{\scriptscriptstyle{G}}(x)-V_0(x)\big)-Z^{\scriptscriptstyle{G_0}}_n(1)\bigg)\mid  Z_1, \ldots, Z_n  \bigg] \bigg| \\
		&   \quad \quad \leq \bigg| \frac{\sqrt{n}}{\sqrt{\log n}}\left(\hat{V}_n(x)-V_0(x)\right)-Z^{\scriptscriptstyle{G_0}}_n(1) \bigg| \xrightarrow{G_0} 0.
	\end{align*} 
	This means that the claim is obtained as soon as we show:
	\begin{align*}
		\sup _{f \in B L_1}\left| \, \mathrm{E}\left(f\bigg(\frac{\sqrt{n}}{\sqrt{\log n}}\big(\hat{V}_{\scriptscriptstyle{G}}(x)-V_0(x)\big)-Z^{\scriptscriptstyle{G_0}}_n(1)\bigg) \mid Z_1,\ldots,Z_n \right)- \ex \left(f(X) \right) \right| \xrightarrow{G_0} 0.	
	\end{align*}
	By Theorem 1.13.1 (iii) in \cite{5} the last display is implied by \eqref{eq: equiv to claim}.
	\end{proof}
	
	\begin{proof}[\hypertarget{proof of reduction}{\textbf{Proof of}} \eqref{eq: representation without remainders}]
		First note that by assumption $Qf$ is a.s.\ well-defined (c.f.\ Remark 4.4 in \cite{4}). Moreover note that \( \ex V_n = |\alpha|/(|\alpha| + n) \) and $|f_n^t(Z)| \leq 2 \sqrt{|t|/ \delta^*_n}$ and \( \forall \, \varepsilon > 0 \), by independence, that:
	\begin{align*}
	\pb \left( \delta_n^{-1} \left| V_n (Q - \mathbb{G}_n) f_t^n \right| > \varepsilon \mid Z_1, \ldots, Z_n \right) &\leq \delta_n^{-1} \varepsilon^{-1} \ex V_n \, \ex \left[  (Q + \mathbb{G}_n) | f_t^n | \mid Z_1, \ldots, Z_n \right] \\
	&\hspace{-0.2cm}\leq 4 \delta_n^{-1} \varepsilon^{-1} \frac{|\alpha|}{|\alpha| + n} \frac{\sqrt{|t|}}{\sqrt{\delta^*_n}}= O\left(\delta^{-1}_n n^{-1} (\delta^{*}_n)^{-1/2} \right).
	\end{align*}
	As $\gamma > \frac{1}{2}$, note that $\delta^{-1}_n n^{-1} (\delta^{*}_n)^{-1/2} = n^{-1} (n / \log{n})^{\frac{1}{2} + \frac{1}{4 \gamma}} \ll  n^{-1} (\sqrt{n}/ \sqrt{\log{n}})^2  \rightarrow 0$. Similarly: $\pb \left( \delta_n^{-1} \left| V_n (\mathbb{B}_n - \mathbb{G}_n) f_t^n \right| > \varepsilon \mid Z_1, \dots, Z_n \right) \stackrel{G_0}{\rightarrow} 0.$

	By Proposition G.2 in \cite{4} for \( \varepsilon_i \overset{\text{i.i.d.}}{\sim} \text{Exp}(1) \) and \( W_{ni} := \varepsilon_i / \bar{\varepsilon}_n \), we have \( (W_{n1}, \ldots, W_{nn}) \sim \text{Dir}(n; 1/n, \dots, 1/n) \). Since \( \mathbb{B}_n \sim \text{DP}(n \mathbb{G}_n) \), for $\delta_{Z_i}$ Dirac measure at $Z_i$, we have the following representation in distribution of $\mathbb{B}_n$ (c.f.\ Example 3.7.9 in \cite{5}):
	\[
	\mathbb{B}_n = \frac{1}{n} \sum_{i=1}^n W_{ni} \delta_{Z_i}.
	\]
	Because of this representation (c.f.\ section 3.7.2. in \cite{5}), $\mathbb{B}_n$ is also known as \textit{Bayesian bootstrap}. Conclude that: 
	\[
		\frac{\sqrt{n}}{\sqrt{\log n}} (\mathbb{B}_n - \mathbb{G}_n) f^n_t = \frac{1}{\sqrt{n \log n}} \sum_{i=1}^n \left( W_{ni} - 1 \right) f^n_t(Z_i).
	\]
	Because $|\bar{\varepsilon}_n - 1| = O_p(1/\sqrt{n})$, for almost every sequence $Z_1,\ldots, Z_n$ conclude that, uniformly in $t$ in compacta: 
	\begin{align}\label{eq: process representation}
		\frac{1}{\sqrt{n \log n}} \sum_{i=1}^n \left( W_{ni} - 1 \right) f_t^n(Z_i)=\frac{1}{\sqrt{n \log n}} \sum_{i=1}^n \left( \varepsilon_i - 1 \right) f_t^n(Z_i) + o_p(1),
	\end{align}
	conclude by recalling the definition of $\mathbb{Z}^*_n$ in \eqref{eq: mathbbZ process}.
	\end{proof}

	\begin{lemma}\label{lemma: cond asymp equicontinuity}
	The process $\mathbb{Z}_n^*$ is conditionally asymptotically equicontinuous, i.e.\ $\forall \, \eta >0$ as $n \rightarrow \infty$ followed by $\delta >0$:
	\begin{align*}
		\pb \bigg\{ \sup _{|s-t|< \delta} \big| \mathbb{Z}_n^*(s) - \mathbb{Z}_n^*(t) \big|>\eta \mid Z_1, \ldots, Z_n\bigg\} \stackrel{G_0}{\longrightarrow} 0.
	\end{align*}
	\end{lemma}
	\begin{proof}
		Take the class of functions for $\delta >0, \: M>0$:  
		\begin{align*}
		\mathscr{F}_{n, \delta}^M :=\big\{ (\log n)^{-1/2} (\varepsilon -1) (f_s^n - f_t^n)(z) : \: s, t \in I_x,|s-t|<\delta, \max (|s|,|t|)<M\big\}.  \numberthis \label{eq: class of functions} 
		\end{align*}
		and the class of functions $\mathscr{H}_{n, \delta}^M$, whose functions for $f \in \mathscr{F}_{n, \delta}^M$ are given by $f/(\varepsilon-1)$, thus without the $\varepsilon-1$. These classes have square integrable natural envelopes $F_{n,\delta}$ and $H_{n,\delta}$ that satisfy for all $\delta>0$ small (cf.\ Lemma 5 in \cite{27}):
		\begin{align*}
			\| F_{n,\delta} \|_2^2 = \ex (\varepsilon-1)^2  \| H_{n,\delta} \|_{2,\scriptscriptstyle{G_0}}^2 =\| H_{n,\delta} \|_{2,\scriptscriptstyle{G_0}}^2  \leq g_0(x + \delta^*_n s \wedge t) \, \delta^2 
		(\log{n})^{-1} ( \log{(\delta \delta^*_n)^{-2}}).
		\end{align*}
		For fixed $s, t \in I_x$, we prove  conditional asymptotic equicontinuity of $\mathbb{Z}_n^*$, i.e.\ $\forall \, \eta >0$ as $n \rightarrow \infty$ followed by $\delta >0$, the right-hand side of the display below converges to 0 under $G_0$:
		\begin{align*}
			\pb \bigg\{ \sup _{|s-t|< \delta} \big| \mathbb{Z}_n^*(s) - \mathbb{Z}_n^*(t) \big|>\eta \mid Z_1, \ldots, Z_n\bigg\} \leq \frac{1}{\eta} \ex \bigg[ \sup _{|s-t|< \delta} \big| \mathbb{Z}_n^*(s) - \mathbb{Z}_n^*(t) \big| \mid Z_1, \ldots, Z_n\bigg].
		\end{align*}
		By the maximal inequality in Theorem 2.14.1 in \cite{5}:
		\begin{align*}
			\ex \bigg[ \sup _{|s-t|< \delta} \big| \mathbb{Z}_n^*(s) - \mathbb{Z}_n^*(t) \big| \bigg] \lesssim J(1) \| F_{n,\delta} \|_2  \lesssim \sqrt{\|g_0\|_{\infty} \, \delta^2 
		(\log{n})^{-1} ( \log{(\delta \delta^*_n)^{-2}})} \to 0,
		\end{align*}
		where $J(1) := J(1,\mathscr{F}_{n, \delta}^M, \mathbb{L}_2 ) = \sup_{Q} \int_{0}^{1} \sqrt{\log N(\eta \| F_{n,\delta} \|_2,\mathscr{F}_{n, \delta}^M, \mathbb{L}_2(Q))} \, d \eta$ and where the square $\mathbb{L}_2(Q)$-distance is given for $s,t$ and $s^\prime, t^\prime \in I_x$ by:
		$
		\frac{1}{\log n} \int ((f^n_s - f^n_t)(z) - (f^n_{s^\prime} - f^n_{t^\prime})(z))^2 \, d Q(z).
		$
		(Note that we did not take $\varepsilon-1$ into account as $\ex (\varepsilon-1)^2=1$).
		Furthermore, we used that $J(1,\mathscr{F}_{n, \delta}^M, \mathbb{L}_2 ) < \infty$ as a consequence of Lemma 2.1 in \cite{27}, where it is shown that the function classes under consideration are VC of uniform bounded index.
	\end{proof}

	In the following Lemma, we make use of the concept of weak convergence of a process in $\ell^{\infty}(K)$ for every compact $K \subset \mathbb{R}$ and $\ell^{\infty}$ being the space of bounded functions on $K$ equipped with the supremum norm. For a review of the conditional weak convergence in $\ell^{\infty}$, see \cite{5} (c.f.\ Example 1.5.1 and Theorem 1.5.4).

	\begin{lemma}[Conditional argmax continuous mapping]\label{lemma: cond argmax}
		Let $\mathbb{Z}_n, \mathbb{Z}$ be stochastic processes indexed by $H \subseteq \mathbb{R}$ such that in probability given a sequence of random variables $Z_1,\ldots,Z_n$:
		\begin{align*}
			\mathbb{Z}_n \mid Z_1,\ldots, Z_n \rightsquigarrow \mathbb{Z} \quad \text{in} \quad \ell^{\infty}(K)
		\end{align*}
		for every compact $K \subset H$. Suppose that almost all sample paths $t \mapsto \mathbb{Z}(t)$ are upper semicontinuous and possess a unique maximum at a (random) point $\hat{t}$, which as a random map in $H$ is tight. If $\, \forall \, n \in \mathbb{N}$, $\mathbb{Z}_n (\hat{t}_n) = \sup_t \mathbb{Z}_n (t)$ a.s.\ and there exist a compact set $\tilde{K} \subset H$ such that, under $\pb$:
			\begin{align*}
				\pb \left( \hat{t}_n \notin \tilde{K} \mid Z_1,\ldots, Z_n  \right) \stackrel{\pb}{\rightarrow} 0,
			\end{align*}
		then $\hat{t}_n \mid Z_1,\ldots,Z_n \rightsquigarrow \hat{t}$ in probability in $H$.
	\end{lemma}
	\begin{proof}
		By the continuous mapping theorem, for every closed set \( F \) and \( A, B \) arbitrary sets of \(H\): $\sup_{t \in F \cap A} \mathbb{Z}_n(t) - \sup_{t \in B} \mathbb{Z}_n(t)$ converges conditionally on \(Z_1, \ldots, Z_n\) in distribution to $\sup_{t \in F \cap A} \mathbb{Z}(t) - \sup_{t \in B} \mathbb{Z}(t),$ where \( \mathbb{Z}(t)\) has a.s. continuous sample paths and a.s. unique maxima.
		Thus using the assumption on $\hat{t}_n$:
		\begin{align*}
		\limsup_{n \to \infty} \pb(\hat{t}_n \in F \cap A \mid Z_1, \ldots, Z_n)
		&\leq \limsup_{n \to \infty} \pb\left( \sup_{t \in F \cap A} \mathbb{Z}_n(t) > \sup_{t \in B} \mathbb{Z}_n(t) \mid Z_1, \ldots, Z_n \right) \\
		&= \pb\left( \sup_{t \in F \cap A} \mathbb{Z}(t) > \sup_{t \in B} \mathbb{Z}(t) \right)
		\leq \pb(\hat{t} \in F \cup B^c),
		\end{align*}
		where in the last inequality we use that \( \mathbb{Z}(\hat{t}) > \sup_{t \notin F^c, t \in A} \mathbb{Z}(t)\) and thus the event 
		\[
		\left\{\sup_{t \in F \cap A} \mathbb{Z}(t) > \sup_{t \in B} \mathbb{Z}(t)\right\} \subseteq \{ \hat{t} \in F \cup B^c\}.
		\]
		Because \(Z_n\) converges in distribution conditionally on \(Z_1, \ldots, Z_n\) to \(Z\) in \(\ell^\infty(K)\) for every compact \(K\), take now \(A = B = K\). Then for every open set \(G\) containing \(\hat{t}\), it holds
		\[
		\mathbb{Z}(\hat{t}) > \sup_{t \notin G, t \in K} \mathbb{Z}(t),
		\]
		because of the continuity of \(\mathbb{Z}\) and the uniqueness of \(\hat{t}\). Thus, the above yields:
		\begin{align*}
		&\limsup_{n \to \infty} \pb(\hat{t}_n \in F \mid Z_1, \ldots, Z_n) 
		\leq \limsup_{n \to \infty} \pb(\hat{t}_n \in F \cap K \mid Z_1, \ldots, Z_n) 
		 \\
		&+ \limsup_{n \to \infty} \pb(\hat{t}_n \in K^c \mid Z_1, \ldots, Z_n) = \pb(\hat{t} \in F \cup K^c).
		\end{align*}
		Note: $\pb(\hat{t} \in F \cup K^c) \leq \pb(\hat{t} \in F) + \pb(\hat{t} \in K^c),$
		where by the stochastic boundedness of \(\hat{t}_n \mid Z_1, \ldots, Z_n\) and the choice of \(K^c\), the terms \(\pb(\hat{t} \in K^c)\) can be made arbitrarily small and $\limsup_{n \to \infty} \pb(\hat{t}_n \in K^c \mid Z_1, \ldots, Z_n)$ converges to $0$.		
	\end{proof}

	\begin{lemma}[Consistency]\label{lemma: consistency}
		Assume \eqref{eq: condition roughness} holds true. The IIP is consistent, i.e. $\forall \, \varepsilon >0$ as $n \rightarrow \infty$:
		\begin{align}\label{eq: consistency}
			\Pi_n (| \hat{V}_{\scriptscriptstyle{G}}(x) - V_0(x)| > \varepsilon \mid Z_1,\ldots, Z_n) \stackrel{G_0}{\rightarrow} 0.
		\end{align}
		Moreover:
		\begin{align}\label{eq: consistency IIP}
			\Pi_n \left( \left| T_{\scriptscriptstyle{G}} \left( V_0(x) + \delta_n  (a + Z^{\scriptscriptstyle{G_0}}_n(1)) \right) - x \right| > \varepsilon \mid Z_1, \ldots, Z_n \right) \stackrel{G_0}{\to} 0.
		\end{align}
	\end{lemma}
	\begin{proof}
		We begin by showing \eqref{eq: consistency}. Let $\psi_x(z) := 2 \sqrt{z} - 2 \sqrt{(z-x)_+}$, using Proposition 4.3 in \cite{4} with $\psi_x$ we obtain:
		\begin{align*}
			&\ex \, (U_{\scriptscriptstyle{G}}(x) \mid Z_1,\ldots,Z_n) = \frac{1}{|\alpha|+n}\int \psi_x  \, d \, \bigg(\alpha + \sum_{i=1}^n \delta_{Z_i} \bigg), \\
			&\mathrm{Var} \, (U_{\scriptscriptstyle{G}}(x) \mid Z_1,\ldots,Z_n) \lesssim \frac{1}{1+|\alpha|+n} \frac{1}{|\alpha|+n} \int \psi_x^2 \,  d \, \bigg(\alpha + \sum_{i=1}^n \delta_{Z_i} \bigg).
		\end{align*}
		Thus, using the conditional Chebyshev inequality:
		\[
		\Pi_n \left( \lvert U_{\scriptscriptstyle{G}} (x) - U_0(x) \rvert > \varepsilon \mid Z_1, \ldots, Z_n \right) 
		\rightarrow 0, \quad \text{a.s.}
		\]
		Note now that the functions \(U_{\scriptscriptstyle{G}}\) are non-decreasing, \(U_0\) is bounded continuous non-decreasing, thus \(\forall \, \varepsilon > 0\) by Exercise 3.2 in \cite{2} we conclude: 
		$$\Pi_n \left( \| U_{\scriptscriptstyle{G}} - U_0 \|_{\infty} > \varepsilon \mid Z_1, \ldots, Z_n \right)  \to 0, \quad \text{a.s.}$$
		Thus, by Marshall's inequality (c.f.\ Exercise 3.1 in \cite{2}), with $ U_{\scriptscriptstyle{G}}^*$ least concave majorant of $ U_{\scriptscriptstyle{G}}$:
		\[
		\Pi_n \left(\lVert U_{\scriptscriptstyle{G}}^* - U_0 \rVert_\infty > \varepsilon \mid Z_1, \ldots, Z_n \right) \stackrel{G_0}{\rightarrow} 0.
		\]
		Now fix $\varepsilon >0$ and $h >0$ such that $| (U_0(x)- U_0(x-h))/h - V_0(x)| < \varepsilon$. By concavity of $U^*_{\scriptscriptstyle{G}}$ we have:
		\begin{align*}
			V_{\scriptscriptstyle{G}}(x) \leq \frac{U^*_{\scriptscriptstyle{G}}(x) - U^*_{\scriptscriptstyle{G}}(x-h)}{h}.
		\end{align*}
		Thus by Lemma 3.1 in \cite{2}:
		\begin{align*}
			&\Pi_n \, (V_{\scriptscriptstyle{G}}(x) \geq V_0(x) + 2 \varepsilon \mid Z_1,\ldots, Z_n) \\
			&\leq \Pi_n \, \left( \frac{U^*_{\scriptscriptstyle{G}} (x) - U^*_{\scriptscriptstyle{G}} (x-h)}{h} \geq \frac{U_{0} (x) - U_{0} (x-h)}{h} +\varepsilon \mid Z_1, \ldots, Z_n \right) \to 0, \quad \text{a.s.} 
		\end{align*}		
		Similarly $V_{\scriptscriptstyle{G}}(x) \geq (U^*_{\scriptscriptstyle{G}}(x+h) - U^*_{\scriptscriptstyle{G}}(x))/h$ gives $\Pi_n \, (V_{\scriptscriptstyle{G}}(x) \leq V_0(x) - 2 \varepsilon \mid Z_1,\ldots, Z_n) \to 0, \quad \text{a.s.}$
		This being true for every $\varepsilon >0$ gives \eqref{eq: consistency}. 
		
		From \eqref{eq: consistency} follows consistency of \(T_{\scriptscriptstyle{G}}\); by the switch relation, \ $\, \forall \, \varepsilon >0$:
		\begin{align*}
		&\Pi_n \left( \left| T_{\scriptscriptstyle{G}} \left( V_0(x) + \delta_n  (a + Z^{\scriptscriptstyle{G_0}}_n(1)) \right) - x \right| > \varepsilon \mid Z_1, \ldots, Z_n \right) \\
		&\leq \Pi_n \left( T_{\scriptscriptstyle{G}} \left( V_0(x) + \delta_n (a + Z^{\scriptscriptstyle{G_0}}_n(1)) \right)   > x+\varepsilon \mid Z_1, \ldots, Z_n \right)  \\
		&\quad + \Pi_n \left(  T_{\scriptscriptstyle{G}} \left( V_0(x) + \delta_n (a + Z^{\scriptscriptstyle{G_0}}_n(1)) \right)  < x - \varepsilon \mid Z_1, \ldots, Z_n \right) \\
		&\leq \Pi_n \left( \hat{V}_{\scriptscriptstyle{G}} (x+\varepsilon) > V_0(x) + \delta_n (a + Z^{\scriptscriptstyle{G_0}}_n(1))  \mid Z_1, \ldots, Z_n \right)  \\
		&\quad + \Pi_n \left( \hat{V}_{\scriptscriptstyle{G}} (x-\varepsilon) < V_0(x) + \delta_n  (a + Z^{\scriptscriptstyle{G_0}}_n(1))   \mid Z_1, \ldots, Z_n \right) \stackrel{G_0}{\to} 0, \numberthis \label{eq: consistency of T_G}
		\end{align*}
		where by \eqref{eq: consistency} above $\hat{V}_{\scriptscriptstyle{G}} (x \pm \varepsilon) \mid Z_1, \ldots, Z_n \rightarrow V_0(x \pm \varepsilon)$, $Z^{\scriptscriptstyle{G_0}}_n(1) \rightsquigarrow \mathbb{Z}(1)$ and $\delta_n \downarrow 0$, thus $\delta_n  (a + Z^{\scriptscriptstyle{G_0}}_n(1)) \stackrel{G_0}{\to} 0$ and $V_0$ is strictly decreasing at $x$ thus we have the strict inequalities $ V_0(x+ \varepsilon) < V_0(x) <  V_0(x- \varepsilon)$.  
	\end{proof}

	\begin{proof}[\hypertarget{proof of stoch boundedness}{\textbf{Proof of}} \eqref{eq: conditional stoch boundedness}]
		The claim is implied by the unconditional version of the convergence given in \eqref{eq: conditional stoch boundedness} and this entails proving that the rate for $\left| T_{\scriptscriptstyle{G}} \left( V_0(x) + \delta_n (a + Z^{\scriptscriptstyle{G_0}}_n(1)) \right) - x \right|$ is $\left( \delta_n^* \right)^{-1}$. Because of Lemma \ref{lemma: consistency}, we can proceed by verifying the rest of the conditions of Theorem 3.2.5 in \cite{5} (or see Theorem 3.4.1 for a more general version). 

		The argmax process we are interested in is $\tilde{T}_n :=  \argmax_t \big\{ \tilde{Z}^{\scriptscriptstyle{G}}_n(t)-a t \big\}$ where $\tilde{Z}^{\scriptscriptstyle{G}}_n$ is given in \eqref{eq: first components Z tilde}-\eqref{eq: second components Z tilde}. Thus note that we can rewrite it as: $\tilde{T}_n = (\delta^*_n)^{-1} \argmax_t \big\{ \mathbb{M}_n(t)\big\} $, where:
		\begin{align}\label{eq: mathbbM_n}
			\mathbb{M}_n(t) := \, \, &2 \int \bigg( \sqrt{(z-x)_+} - \sqrt{(z-x -t)_+} - \frac{1}{2} V_0(x)t \bigg) \, d (G - \mathbb{G}_n+G_0)(z)\\
			&- 2 \frac{t}{\delta^*_n}  \int \bigg( \sqrt{(z-x)_+} - \sqrt{(z-x -\delta^*_n)_+} \bigg) \, d (\mathbb{G}_n-G_0)(z) \\
			&+ 2  \int \bigg( \sqrt{(z-x)_+} - \sqrt{(z-x -t)_+} \bigg) \, d (\mathbb{G}_n-G_0)(z)  - at \delta_n .
		\end{align}
		Moreover define the deterministic function:
		\begin{align*}
			\mathbb{M}(t):= 2 \int \bigg( \sqrt{(z-x)_+} - \sqrt{(z-x -t)_+} - \frac{1}{2} V_0(x)t \bigg) \, d G_0(z).
		\end{align*}
		Note that $\mathbb{M}_n(0) = \mathbb{M}(0) = 0$. In \cite{27} (c.f.\ proof of (18)) it is shown that $\exists \, \delta_0>0$ such that $\forall \, |t|< \delta_0$:
		\begin{align}\label{eq: bound curvature}
			\mathbb{M}(t) - \mathbb{M}(0) \leq - K \, |t|^{\gamma +1}
		\end{align} 
		Following Theorem 3.2.5 in \cite{5}, define thus the metric $d(s,t) := |s-t|^{\frac{\gamma+1}{2}}.$

		Consider the class of functions $h_t :[0, \infty) \rightarrow \mathbb{R}$:
		\begin{align*}
		\mathscr{H} := \big\{h_t \: : \: h_t (z) =& \: 2 \sqrt{(z-x)_{+}} -2 \sqrt{(z-x-t)_{+}}  - V_0(x)t, \: t \geq 0 \big\}
		\end{align*}
		and, for $\eta>0$, its subclass $ \mathscr{H}_\eta :=\left\{h_t \in \mathscr{H}:|t| \leq \eta\right\}$. In \cite{1} it is shown that the envelope of $\mathscr{H}_\eta$, $H_{\eta} (z) := 2 \sqrt{(z-x)_{+}} -2 \sqrt{(z-x-\eta)_{+}} - V_0(x)\eta$, satisfies for all $\eta$ small enough (c.f.\ (11) in \cite{1}):
		\begin{align}\label{eq: bound envelope 1}
			\int H^2_{\eta}(z) \, g_0(z) \, dz \leq -2g_0(x) \eta^2 \log{(\eta)}.
		\end{align}

		Consider now the class of functions $\widetilde{\mathscr{H}}_n := \{ \tilde{h}^n_t \, : \, t \geq 0 \}$ where, for any fixed $t \geq 0$:
		\begin{align*}
		\tilde{h}^n_t (z) := \:  2 \bigg(\frac{t}{\delta^*_n} -1\bigg) \sqrt{(z-x)_{+}} -2\frac{t}{\delta^*_n} \sqrt{(z-x-\delta^*_n)_{+}}   +2 \sqrt{(z-x-t)_{+}} ,
		\end{align*}
		and, for $\eta>0$, its subclass $ \widetilde{\mathscr{H}}_{n,\eta}:=\big\{\tilde{h}^n_t \in \widetilde{\mathscr{H}}_n:|t| \leq \eta\big\}$. Using the same computations as in \cite{1}, it follows that the envelope of $\widetilde{\mathscr{H}}_{n,\eta}$, $\tilde{H}_{n,\eta} (z) := \sup_{|t|< \eta }|\tilde{h}^n_t(z)|$, satisfies for all $\eta$ small enough:
		\begin{align}\label{eq: bound envelope 2}
			\int \tilde{H}^2_{n,\eta} (z) \, g_0(z) \, dz \leq -2g_0(x) \eta^2 (\log{(\eta)} + \log{(\delta^*_n)}).
		\end{align}

		Using the representation \eqref{eq: representation process G} we have:
		\begin{align*}
			\mathbb{M}_n(t) - \mathbb{M}(t) &=  (1-V_n) \int h_t(z) \, d (\mathbb{B}_n - \mathbb{G}_n)(z)   + V_n \int h_t(z) \, d (Q- \mathbb{G}_n)(z) - at \delta_n \\
			& \quad -  \int \tilde{h}^n_t(z) \, d (\mathbb{G}_n-G_0)(z), 
		\end{align*}
		and, by applying the maximal inequality in Theorem 2.14.1 in \cite{5} for sufficiently small $\delta >0$,  we conclude ($d(t,0)< \delta \, \iff \, |t|< \delta^{\frac{2}{1+\gamma}}$):
		\begin{align*}
			&\ex \bigg[ \sup_{|t|< \delta^{\frac{2}{1+\gamma}}} \sqrt{n} |(\mathbb{M}_n - \mathbb{M})(t)| \bigg] \lesssim J(1, \mathscr{H}_{\delta^{\frac{2}{1+\gamma}}}, \mathbb{L}_2) \sqrt{\int H^2_{\delta^{\frac{2}{1+\gamma}}} \, dG_0} + \sqrt{n} \frac{1}{n} \delta^{\frac{2}{1+\gamma}} \\
			&\quad + |a| \delta^{\frac{2}{1+\gamma}} \delta_n \sqrt{n} + J(1, \widetilde{\mathscr{H}}_{n,\delta^{\frac{2}{1+\gamma}}}, \mathbb{L}_2) \sqrt{\int \tilde{H}^2_{n,\scriptscriptstyle{\delta^{\frac{2}{1+\gamma}}}}\, dG_0 \,}  \lesssim  \delta^{\frac{2}{1+\gamma}} \big(\sqrt{\log{\delta}} + \sqrt{\log{n}} \big).
		\end{align*}
		The first term after applying the maximal inequality comes from the $(1-V_n)(\mathbb{B}_n - \mathbb{G}_n)$-term in the representation above of $\mathbb{M}_n - \mathbb{M}$; the second one comes from the $V_n(Q - \mathbb{G}_n)$-term, the third one from the $at \delta_n $-term and the last one from the $\mathbb{G}_n-G_0$-term (having used also for this one the maximal inequality). The last inequality is just a consequence of \eqref{eq: bound envelope 1}, \eqref{eq: bound envelope 2} and the facts that $J(1, \mathscr{H}_{\delta^{\frac{2}{1+\gamma}}}, \mathbb{L}_2) < \infty$ by the proof of Lemma 2.1 in \cite{27} and by the same reasonings $J(1, \widetilde{\mathscr{H}}_{n,\delta^{\frac{2}{1+\gamma}}}, \mathbb{L}_2) $ is uniformly bounded.

		The above computations and Theorem 3.2.5 in \cite{5} imply that the rate for $d$ is determined as follows: $\delta^{\frac{2}{1+\gamma}} (\sqrt{\log{\delta}} + \sqrt{\log{n}}) \leq \sqrt{n} \delta^2 \ $, hence $ \delta \gtrsim \left( \frac{\log{n}}{n}\right)^{\frac{1+\gamma}{4 \gamma}}$.
		Because $d(s,t) = |s-t|^{\frac{\gamma+1}{2}}$, the rate for the distance $|\cdot|$ is given by: 
		$$\left( \frac{\log{n}}{n}\right)^{\frac{1+\gamma}{4 \gamma} \cdot \frac{2}{1+\gamma}} = \left( \frac{\log{n}}{n}\right)^{\frac{1}{2\gamma}} \equiv \delta^*_n.$$
		Altogether the above proves that the rate for $\argmax_t \big\{ \mathbb{M}_n(t)\big\}$ is $(\delta^*_n)^{-1}$.
	\end{proof}

\end{appendix}
%
%

\begin{acks}[Acknowledgments]
We thank the anonymous reviewers and the Associate Editor for their thoughtful comments and constructive suggestions, which have helped us improve the clarity and quality of this work.
\end{acks}
%

\begin{supplement}
\stitle{Supplement to the article: ``Semiparametric Bernstein-von Mises Phenomenon via Isotonized Posterior in Wicksell's problem'' \textnormal{(available at \cite{64})}.}
\end{supplement}


\bibliographystyle{imsart-number} 
\bibliography{mybiblio}       


\end{document}